\documentclass[a4paper]{article}
\usepackage{amssymb,amsbsy,amsmath,amsfonts,amssymb,amscd}
\usepackage{latexsym}
\usepackage{geometry}
\usepackage{latexsym}

\newcommand{\lie}[1] {\mathfrak{#1}}  
\newcommand{\bb}[1]{{\mathbb #1}}    

\newcommand{\GL}{{\rm GL}}
\newcommand{\SL}{{\rm SL}}

\newcommand{\Aut}{{\rm Aut}}
 \newcommand{\PSL}{{\rm PSL}}

\newcommand{\A}{{\rm A}}

{ \par  \medskip  \noindent  {\bf Definition} \hspace{0.5em} }%
{\par \medskip } 

\newtheorem{example1}{Example}[section]

{ \begin{example1} \rm }%
{ \end{example1}} 
\newenvironment {proof}%
{ \noindent {\em Proof: }}%
{\hspace*{\fill}$\Box$\par \medskip } 
{ \noindent {\em #1 \hspace{0.5em}} }%
{\hspace*{\fill}$\Box$\par \medskip } 
{{\em Remark \hspace{0.5em}} }%
{\par \medskip } 
\newenvironment {prof}%
{ \noindent {\em Proof }}%
{\hspace*{\fill}$\Box$\par \medskip }

\newtheorem{proposition}{Proposition}[section]
\newtheorem{definition1}[proposition]{Definition}
\newtheorem{theorem}[proposition]{Theorem}
\newtheorem{lemma}[proposition]{Lemma}
\newtheorem{corollary}[proposition]{Corollary}
{ \begin{definition1} \rm }%
{ \end{definition1}}

\input xy
\xyoption{all}

\font\Bbb=msbm10

\def\BZ{{\hbox{\Bbb Z}}}

\def\BN{\hbox{\Bbb N}}

\def\BC{{\hbox{\Bbb C}}}
\def\BQ{{\hbox{\Bbb Q}}}

\def\A{{\hbox{\rm Aut}}}
\def\Ro{{\bar R_{\rm old}}}
\def\A{{\rm Aut}}     

\def\op{{\rm op}}

\def \Box{\lower .1 em
          \vbox{\hrule \hbox{\vrule \hskip .6 em \vrule height .6 em} \hrule}}
\def \Mid{\quad \vrule \quad}
\font \msbm= msbm10
\def \rfish{\mathbin {\hbox {\msbm \char'157}}}

\title{ Linear Representations of the Automorphism Group of a 
Free Group}

\author{Fritz Grunewald\\
Mathematisches Institut \\ 
Heinrich-Heine-Universit\"at\\ 
D-40225 D\"usseldorf\\
e-mail: fritz@math.uni-duesseldorf.de 
\and 
Alexander Lubotzky\\
Einstein Institute of Mathematics\\
The Hebrew University of Jerusalem\\
Jerusalem, 91904, Israel\\
e-mail: alexlub@math.huji.ac.il}

\begin{document}
\maketitle

\begin{abstract} 
Let $F_n$ be the free group on $n\ge 2$ elements and $\A(F_n)$ its group of
automorphisms. In this paper we present a rich collection of 
linear representations of
$\A(F_n)$ arising through the action of finite index subgroups of it on
relation modules of finite quotient groups of $F_n$. We show 
(under certain conditions) that the images
of our representations are arithmetic groups. 
\newline 

\end{abstract}

\medskip

\bigskip
\bigskip

\tableofcontents 
\medskip

\bigskip

\section{Introduction} 

Let $F_n$ be the free group on $n\ge 2$ elements and $\A(F_n)$ its group of
automorphisms. The latter is a much studied group, but very little seems to be
known about its (finite dimensional) complex representation theory
(see for example \cite{FP}, \cite{DJ}, \cite{PORA}, \cite{RA}, \cite{LP}). 
In fact,
as far as we know the only representations studied (at least when $n\ge 3$) 
are:
\begin{itemize}
\item The representation
\begin{equation}\label{furep}
\rho_1 : \A(F_n)\to \A(F_n/F'_n)\cong \GL(n,\BZ) 
\end{equation}
where $F_n'$ is the commutator subgroup of $F_n$.
\item The representations factoring through the homomorphisms
$$\rho_i : \A(F_n)\to \A(F_n/F^{(i+1)}_n)$$
where $F^{(i+1)}_n$ ($i=0,1,\ldots$) stands for the lower 
central series of $F_n$.  
\item Representations with finite image. 
\end{itemize}
What is common to all of the above mentioned reprentations $\rho$ is that if we
denote by ${\cal H}:=\overline{\rho(\A(F_n))}$ the Zariski closure of 
$\rho(\A(F_n))$, then the semi-simple part of the connected component 
${\cal H}^\circ$ of ${\cal H}$ is either trivial or $\SL(n,\BC)$. Moreover,
$\rho({\rm IA}(F_n))$ is always virtually solvable for these representations,
where ${\rm IA}(F_n)$ is the kernel of $\rho_1$.

In this paper we will show that the representation theory of $\A(F_n)$ is in
fact much richer. Our main Theorem (Theorem \ref{teo} below) 
implies for example 
\begin{theorem}\label{teo4}
Let $n\ge 2$, $k\ge 1$, $h_1<\ldots <h_k$, $m_1,\ldots ,m_k$ be 
natural numbers.  Let $\BQ(\zeta_{m_i})$ be the field of 
$m_i$-th roots of unity and
$\BZ(\zeta_{m_i})$ its ring of integers. There is a 
subgroup $\Gamma\le \A(F_n)$ of finite index and a representation 
$$
\rho : \Gamma \to \prod_{i=1}^k\SL((n-1)h_i,\BQ(\zeta_{m_i}))^{m_i} 
$$
such that $\rho(\Gamma)$ is commensurable with 
$\prod_{i=1}^k\SL((n-1)h_i,\BZ(\zeta_{m_i}))^{m_i}$. 
\end{theorem}
Theorem \ref{teo4} shows that $\A(F_n)$ has a representation $\rho$ such that
the connected component of $\overline{\rho(\A(F_n))}$ is isomorphic to 
$$\prod_{i=1}^k\SL((n-1)h_i,\BC)^{m_i}$$
as above. In fact the proof shows that this is even true for 
$\overline{\rho({\rm IA}(F_n))}$ and in particular the latter 
is very far from being virtually solvable (see Section \ref{Concl1} for more). 

Specifying $m_1=\ldots =m_k=1$ in Theorem \ref{teo4}, we deduce
\begin{corollary}\label{coro4}
Let $n\ge 2$, $k\ge 1$ be natural numbers.  There is a 
subgroup $\Gamma\le \A(F_n)$ of finite index and a representation 
$$
\rho : \Gamma \to \prod_{i=1}^k\SL((n-1)i,\BZ) 
$$
such that $\rho(\Gamma)$ is of finite index in $\prod_{i=1}^k\SL((n-1)i,\BZ)$. 
\end{corollary}
Specialising Theorem \ref{teo4} even further to 
$n=3$, $k=1$, $h_1=m_1=1$ we get that
$\A(F_3)$ has a subgroup of finite index which can be mapped onto a subgroup
of finite index in $\SL(2,\BZ)$. This implies  
\begin{corollary}\label{kaz}
The automorphism group $\A(F_3)$ is large, that is it has a subgroup of finite
index which can be mapped onto a free nonabelian group. In particular
$\A(F_3)$ does not have Kazdhan's property (T).
\end{corollary}
The corollary is easy for $\A(F_2)$, see the discussion in Section 
\ref{Proofs}, but it is a
well known open problem for larger $n$. Our solution for 
$n=3$ does not indicate what should be the answer for $n>3$. 

We shall describe now how the representations to be considered here arise.  
Let $G$ be a finite group and $\pi : F_n\to G$ a surjective homomorphism
of the free group $F_n$ onto $G$. Let $R$ be the kernel of 
$\pi$. We define  
\begin{equation}
\Gamma(R):=\{\, \varphi\in \A(F_n)\Mid \varphi(R)=R \, \}
\end{equation}
and 
\begin{equation}\label{I1}
\Gamma(G,\pi):=\{\, \varphi\in \Gamma(R) \Mid  
\varphi\ {\rm induces\ the\ identity\ on}\ F_n/R \, \}.
\end{equation}
These are subgroups of finite index in $\A(F_n)$. We let
\begin{equation}
\bar R:=R/R'
\end{equation}
be the relation module of the presentation $\pi : F_n\to G$. 
The action of $F_n$ on $R$ by conjugation leads 
to an action of $G$ on $\bar R$.

The structure of $\bar R$ as a $G$-module is described by 
Gasch\"utz' theory (see Section \ref{suga} and the references therein). 
It says in particular that 
\begin{equation}\label{G1}
\BQ\otimes_{\bb Z} \bar R\cong \BQ\oplus \BQ[G]^{n-1}
\end{equation}
as a module over the rational group ring $\BQ[G]$ (see also (\ref{gasch})).
Our work can be described as an equivariant Gasch\"utz' theory where we try to
describe the action of $\Gamma(G,\pi)$ on $\bar R$: Every $\varphi\in
\Gamma(R)$ induces a linear automorphism $\bar \varphi$ 
of $\bar R$ and (as an application of the theorem of Gasch\"utz) 
$\Gamma(G,\pi)$ consists exactly of those elements  $\varphi\in\Gamma(R)$ 
for which $\bar\varphi$ is $G$-equivariant.
Already now we can see from (\ref{G1}) why the number $n-1$ plays such an
important role in Theorem \ref{teo4}.

To be more precise: 
The relation module $\bar R$ is a finitely generated 
free abelian group, let $t$ be its $\BZ$-rank. Then
\begin{equation}\label{I4}
{\cal G}_{G,\pi}:=\Aut_G(\BC\otimes_{\bb Z} \bar R)\le \GL(t,\BC)
\end{equation}
is a $\BQ$-defined algebraic subgroup of $\GL(t,\BC)$, in fact 
${\cal G}_{G,\pi}$ is the centraliser of the group $G$ acting on 
$\BC\otimes_\BZ \bar R$ through matrices with rational entries. Let
\begin{equation}
{\cal G}_{G,\pi}^1 \le \SL(t,\BC)
\end{equation}
be the kernel of all $\BQ$-defined homomorphisms from the complex algebraic
group $\BC\otimes{\cal G}_{G,\pi}$ to the multiplicative group. This is a
$\BQ$-defined subgroup of ${\cal G}_{G,\pi}$. We shall describe it in more
detail in Section \ref{suga}.
We define 
\begin{equation}
{\cal G}_{G,\pi}(\BZ):=\{\, \phi\in {\cal G}_{G,\pi}\Mid \phi(\bar
R)=\bar R\,\}.
\end{equation}
This is an arithmetic subgroup of the $\BQ$-defined algebraic group
${\cal G}_{G,\pi}$ which contains the arithmetic subgroup
\begin{equation}
{\cal G}^1_{G,\pi}(\BZ):=\{\, \phi\in {\cal G}^1_{G,\pi}\Mid \phi(\bar
R)=\bar R\,\}
\end{equation}
of ${\cal G}^1_{G,\pi}$.
Following the definitions we obtain an 
integral linear representation
\begin{equation}\label{I5}  
\rho_{G,\pi} : \Gamma(G,\pi)\to {\cal G}_{G,\pi}(\BZ),\qquad
\rho_{G,\pi}(\varphi)=\bar\varphi\quad (\varphi\in \Gamma(G,\pi)).
\end{equation}
Examples show (see Section \ref{Choose}) that in general $\rho_{G,\pi}$ is not
onto, neither is its image of finite index in ${\cal G}_{R,\pi}(\BZ)$. 
However, our main result shows that, at least for redundant presentations, the
image of $\rho_{G,\pi}$ captures the whole semi-simple part of 
${\cal G}_{R,\pi}(\BZ)$. More precisely:
\begin{theorem}\label{teo} Assume $n$ is a natural number with $n\ge 4$.
Let $\pi : F_n\to G$ be a redundant presentation of the finite group $G$. Then
$\rho_{G,\pi}(\Gamma(G,\pi))\cap {\cal G}^1_{G,\pi}$ is 
of finite index in the arithmetic group ${\cal G}_{G,\pi}^1(\BZ)$.
\end{theorem}
A presentation $\pi : F_n\to G$ is called {\it redundant} if there is a basis
$x_1,\ldots ,x_n$ of $F_n$ such that $\pi(x_n)=1$. This is equivalent to the
kernel of $\pi$ containing at least one primitive element 
(that is a member of a basis) of $F_n$.

In general we do not know whether 
$\rho_{G,\pi}(\Gamma(G,\pi))\cap {\cal G}_{G,\pi}^1(\BZ)$ is of finite index
in $\rho_{G,\pi}(\Gamma(G,\pi))$. Section \ref{SL} contains criteria under
which this is going to happen.

Choosing the finite group in various ways we obtain a rich variety of 
linear representations of subgroups of finite index in $\A(F_n)$ with
arithmetic image groups. Let us start with some simple examples.

{\it Example 1:} Let $G=C_2$ be the group of order $2$ and $\pi :F_n\to G$
be any surjective homomorphism. Since we generally have assumed $n\ge 2$, the
representation is redundant. As an abelian group we have 
$\bar R\cong \BZ^{2(n-1)+1}$. Also $\BQ\otimes_{\bb Z} \bar R$ 
can be decomposed into
the $+1$-eigenspace of a generator of $G$ and the corresponding
$-1$-eigenspace. We obtain a decomposition 
$\BQ\otimes_{\bb Z} \bar R=\BQ^n\oplus \BQ^{n-1}$ 
which has to be respected by   
$\rho_{G,\pi}(\Gamma(G,\pi))$. Thus we obtain a representation
$$
\rho_{G,\pi} : \Gamma(G,\pi) \to \GL(n,\BQ)\times \GL(n-1,\BQ)
$$
which we show to have an image commensurable with $\SL(n,\BZ)\times
\SL(n-1,\BZ)$.

The first factor is the familiar one arising from the representation
(\ref{furep}). But the second is somewhat less expected. It already gives rise
to a representation which is suitable for proving Corollary \ref{kaz}.

{\it Example 2:} Let $p$ be a prime number, $\zeta_p\in\BC$ a non trivial
$p$-th root of unity. Let $K:=\BQ(\zeta_p)$ be the $p$-th cyclotomic field
and $\BZ(\zeta_p)$ its ring of integers. 
Taking $G=C_p$ the cyclic group of order 
$p$ we have $\bar R\cong \BZ^{p(n-1)+1}$. Now $p(n-1)+1=(p-1)(n-1)+n$ and we  
get a representation
$$
\rho_{G,\pi} : \Gamma(G,\pi) \to \GL(n,\BQ)\times \GL(n-1,K)
$$
which we show to have an image commensurable with  $SL(n,\BZ)\times
SL(n-1,\BZ(\zeta_p))$. This example is treated in Section \ref{Choose} along
with the case of a general cyclic group.

Let $G$ be now an arbitrary finite group and $\pi: F_n\to G$ a redundant
epimorphism. Decompose the group algebra $\BQ[G]$ 
as
$$\BQ[G]=\BQ\oplus \prod_{i=2}^\ell M(h_i,D_i)$$
where the $D_i$ are division algebras. Theorem \ref{teo} implies that there is 
a homomorphism $\rho$ from a finite index subgroup $\Gamma$ of $\A(F_n)$ into
a subgroup of $\prod_{i=2}^\ell\GL((n-1)h_i,D_i)$ where the image contains an
arithmetic subgroup of 
$$\Lambda:=\prod_{i=2}^\ell\SL((n-1)h_i,{\cal R}_i)$$ where 
${\cal R}_i$ is an order in $D_i$. In general we do not know whether
$\rho(\Gamma)\cap\Lambda$ is of finite index in $\rho(\Gamma)$. It will be
especially interesting if it is not! In this case  
$\rho(\Gamma)/(\rho(\Gamma)\cap\Lambda)$ would be an infinite abelian group and
this would imply that $\A(F_n)$ does not have Kazdhan property (T). We have
some results in the opposite direction showing, for example, that if $G$ is
metabelian then $\rho(\Gamma)$ is in fact commensurable to $\Lambda$. 

In any event, whenever $\rho(\Gamma)\cap\Lambda$ is of finite index or not, it
becomes of finite index after we divide by the center of 
$\prod_{i=2}^\ell\GL((n-1)h_i,D_i)$. This shows that in all cases an
arithmetic group commensurable to $\Lambda$ is a quotient of some finite index
subgroup of $\A(F_n)$. Thus Theorem \ref{teo} provides a very rich class of
arithmetic subgroups which are quotients of (finite index subgroups) of
$\A(F_n)$. Theorem \ref{teo4} is obtained from the main result by some special
choices of the finite group $G$ and by an application of Margulis
super-rigidity (see Section \ref{Proofs}).

The present paper opens up a large number of interesting problems. Among
others it shows that the classical Torelli subgroup (i.e. ${\rm IA}(F_n)$) 
is just a first example in a series of infinitely many. 
The kernel of $\rho_{G,\pi}$, for any (redundant) homomorphism 
$\pi:F_n\to G$ of $F_n$ to a finite group $G$, 
can be considered as a natural generalisation of
${\rm IA}(F_n)$ which corresponds to the case where $G$ is the trivial group.
Just as ${\rm IA}(F_n)$, the kernel of $\rho_{G,\pi}$ is residually
torsion-free nilpotent and, at least when $\pi$ is redundant,
the image of $\rho_{G,\pi}$ is an arithmetic group hence
is finitely presented, which implies that   
${\rm ker}(\rho_{G,\pi})$ is a finitely generated normal subgroup.  
Is ${\rm ker}(\rho_{G,\pi})$ a finitely generated group? In case of the
classical Torelli group this is a result of Magnus (see \cite{MKS}).

Let us now describe the method of the proof of the main result and the layout
of the paper: In Section \ref{Two}, we describe in detail the rational
relation module $\BQ\otimes_{\bb Z}\bar R$ together 
with ${\cal G}_{G,\pi}$ which
is the centraliser of $G$ when acting on the relation module. In Section 
\ref{Drei}, we describe the homomorphism from $\Gamma (G,\pi)$ to 
${\cal G}_{G,\pi}$, which is given by the action of $\Gamma (G,\pi)$
on $\BQ\otimes_\BZ\bar R$. We compute the action of various elements of    
$\Gamma (G,\pi)$ on $\BQ\otimes_{\bb Z}\bar R$ in Section \ref{Vier} 
using results
on Fox calculus from Section \ref{Drei}. A careful choice of elements from 
$\Gamma (G,\pi)$ allows us to produce sufficiently many unipotent elements
of ${\cal G}_{G,\pi}(\BZ)$ to appeal to results of Vaserstein \cite{Vas} 
(see also Venkataramana \cite{Venk} for generalisations) to deduce 
in Section \ref{Fuenf} that these
elements generate a finite index subgroup. Section \ref{Choose} is devoted to
a detailed study of the case when $G$ is a cyclic group. In this case we give
fairly complete results. Some of them are used in Section \ref{Proofs} for the
proof of Theorem \ref{teo4}. 
In Section \ref{SL} we discuss the question
whether the image of $\Gamma (G,\pi)$ is commensurable with 
${\cal G}_{G,\pi}^1(\BZ)$ or an infinite extension of it. We solve it
for cyclic groups and also for metabelian groups. 
Some remarks, computational
results and open problems are given in Section \ref{Concl}.

In a subsequent paper \cite{GLu} we apply similar methods to study linear
representations of the mapping class group.

{\it Acknowledgement:} The authors are grateful to A. Rapinchuk for a useful
conversation which led to the results about abelian groups in Section
\ref{SL}. They thank the Minerva-Landau Center for its support and especially the 
Institute for Advanced Study in Princeton, where this
work was done, for its hospitality in the academic year 2005/2006. The second
author thanks the NSF, the Monell Foundation 
and the Ellentuck Fund for support. 

\section{Relation modules and representations of subgroups of finite index in
$\A(F_n)$ }\label{Two}

In this section we collect some results on the structure of the relation
modules of finite groups.  We also consider 
integrality properties of the linear reprentations of
subgroups of finite index in $\A(F_n)$ introduced in the introduction.

\subsection{The rational group ring}\label{sura}

Here $G$ is a finite group and $\BQ[G]$ is its rational group ring. 
Our basic reference for the structural results we need is the book of Curtis
and Reiner \cite{CR} (see also \cite{S}), in particular 
Sections IV and XI therein.  

We start to set up some notation. 
Let $N$ be an irreducible left $\BQ[G]$-module. We write
$$
{\bf I}_N(M)
$$
for the $N$-isotypic component of the left $\BQ[G]$-module $M$. We also choose
representatives $N_1,\ldots , N_\ell$ for the isomorphism classes of 
irreducible left $\BQ[G]$-modules and abreviate 
$$
{\bf I}_i(M):={\bf I}_{N_i}(M)\qquad (i=1,\ldots,\ell)
$$ 
for a left $\BQ[G]$-module $M$. We take $N_1$ to be the one dimensional
trivial module.

Viewing $\BQ[G]$ as a left $\BQ[G]$-module we obtain a direct product
decomposition 
\begin{equation}\label{ri}
\BQ[G]=\prod_{i=1}^\ell\, {\bf I}_i(\BQ[G]).
\end{equation}
In this case every ${\bf I}_i(\BQ[G])$ is also a right submodule of 
$\BQ[G]$. Furthermore every ${\bf I}_i(\BQ[G])$ is an ideal of $\BQ[G]$ having
its own unity.
The decomposition
(\ref{ri}) is a decomposition of $\BQ[G]$ as a direct product of 
$\BQ$-algebras.

We let $P_i\in \BQ[G]$ ($i=1,\ldots,\ell$) be the right projector of 
$\BQ[G]$ onto ${\bf I}_i(\BQ[G])$, that is we have
\begin{itemize}
\item[{\bf RP}:]
$AP_i=A$ for all $A\in {\bf I}_i(\BQ[G])$       
and $AP_i=0$
if $A$ is in the direct product of the subrings
${\bf I}_j(\BQ[G])$ with ${\bf I}_i(\BQ[G])$ excluded.
\end{itemize} 
Notice that we also have 
\begin{itemize}
\item[{\bf LP}:]
$P_iA=A$ for all $A\in {\bf I}_i(\BQ[G])$ and $P_iA=0$ 
if $A$ is in the direct product of the subrings
${\bf I}_j(\BQ[G])$ with ${\bf I}_i(\BQ[G])$ excluded.
\end{itemize}
In fact we have $P_i\in {\bf I}_i(\BQ[G])$ for $i=1,\ldots ,\ell$
and $P_i$ can be thought of as the
unit element in ${\bf I}_i(\BQ[G])$.

Given two left $\BQ[G]$-modules $M_1$, $M_2$ we define 
${\rm Hom}_G(M_1,M_2)$ to be the vector space of $\BQ[G]$-module homomorphisms
from $M_1$ to $M_2$. The $\BQ$-algebra of $\BQ[G]$-module endomorphisms of 
$M_1$ is denoted by ${\rm End}_G(M_1)$.

Let us write $S^{\rm op}$ for the {\it opposite ring} 
of an associative ring $S$. 
Mapping an element  $A\in {\bf I}_i(\BQ[G])$ to the 
multiplication from the right by $A$ defines $\BQ$-algebra isomorphisms
\begin{equation}\label{isoppp}
{\bf I}_i(\BQ[G])^{\rm op}\to 
{\rm End}_G({\bf I}_i(\BQ[G])) \qquad (i=1,\ldots ,\ell).
\end{equation}
The existence of the isomorphisms (\ref{isoppp}) implies that
\begin{equation}\label{dim17}
{\rm dim}_{\bb Q}({\rm Hom}_G({\bf I}_i(\BQ[G])^f,{\bf I}_i(\BQ[G]))
=f{\rm dim}_{\bb Q}({\bf I}_i(\BQ[G]))
\end{equation}
for all $f\in\BN$ and $i=1,\ldots,\ell$.

The algebras 
$$
D_i:={\rm End}_G(N_i)\qquad (i=1,\ldots,\ell)
$$
of left $\BQ[G]$-module endomorphism of the irreducible modules $N_i$ are 
finite dimensional division algebras with cyclotomic number fields as center.
We set 
$$h_i:={\rm dim}_{D_i}(N_i)\qquad (i=1,\ldots,\ell).$$ 
Then the algebra ${\bf I}_i(\BQ[G])$
is isomorphic to a $h_i\times h_i$ matrix algebra $M(h_i,D_i)$ over $D_i$,
that is
\begin{equation}\label{ra1}
{\bf I}_i(\BQ[G])\cong M(h_i,D_i)\qquad (i=1,\ldots ,\ell).
\end{equation}
Using (\ref{isoppp}) we get isomorphisms of algebras
\begin{equation}\label{ra1oppp}
 {\rm End}_G({\bf I}_i(\BQ[G])) \cong M(h_i,D_i^{\rm op})
\qquad (i=1,\ldots ,\ell).
\end{equation}
We shall further discuss integral versions of the isomorphisms (\ref{ri}) and
(\ref{ra1}). Our reference here is \cite{CR} Section XI or \cite{D} 
Section XX. 

Let $K$ be a number field with ring of integers ${\cal O}$ 
and $S$ a finite dimensional 
$K$-algebra. A subring ${\cal R}$ of $S$ is called an {\it order} if 
it contains a $K$-basis of $S$ and if there is a subring ${\cal O}_0$ of
finite index in ${\cal O}$ such that ${\cal R}$ is a finitely generated
${\cal O}_0$-module
(subrings are always assumed to contain the identity element of the bigger
ring). Let $h$ be a natural number and $S=M(h,D)$ where $D$ is a 
finite dimensional division algebra over $K$ with center $K$. We write
${\rm End}_S(S^m)$ ($m\in\BN$) for the algebra of $S$-left module
endomorphisms of the $S$-left module $S^m$. We obtain an algebra isomorphism 
\begin{equation}\label{ord}
\Theta: {\rm End}_S(S^m)\to M(m,S)^\op=M(m,S^\op)=M(m,M(h,D^\op)). 
\end{equation}
Here we identify $M(m,S)^\op$ in the usual way with $M(m,S^\op)$ and then with
$M(m,M(h,D^\op))$.
We call a subgroup $\Lambda$ of the additive group of $S^m$ a {\it lattice} if 
$\Lambda$ contains a $K$-basis of $S^m$ and if 
there is a subring ${\cal O}_0$ of finite index in ${\cal O}$ such that
$\Lambda$ is a finitely generated ${\cal O}_0$-module.
We write
${\rm End}_S(S^m)_{\Lambda}$ for the ring of those endomorphisms stabilising
$\Lambda$. The following is well known (see \cite{CR}, \cite{D}).
\begin{lemma}\label{repkon} Let $\Lambda\le S^m$ ($m\in \BN$) be a lattice.
There is a sublattice $\tilde\Lambda\le \Lambda$ (of finite index in
 $\Lambda$) and an order
 ${\cal R}$ in $D$ such that
$$\Theta({\rm End}_S(S^m)_{\tilde\Lambda})\subset M(m,M(h,{\cal R}^\op)). $$
\end{lemma}
Let ${\rm Aut}_S(S^m)_\Lambda$ be the group of invertible elements in 
${\rm End}_S(S^m)_\Lambda$. Note that 
${\rm Aut}_S(S^m)_{\tilde\Lambda}$ is commensurable with 
${\rm Aut}_S(S^m)_{\Lambda}$ in particular 
${\rm Aut}_S(S^m)_{\Lambda}\cap {\rm Aut}_S(S^m)_{\tilde\Lambda}$
has finite index in ${\rm Aut}_S(S^m)_{\Lambda}$ and we obtain a
representation 
\begin{equation}\label{ord1}
\Theta: {\rm Aut}_S(S^m)_{\Lambda}\cap {\rm Aut}_S(S^m)_{\tilde\Lambda}
\to \GL(m,M(h,{\cal R}^\op))=\GL(mh,{\cal R}^\op).
\end{equation}

\subsection{A result of Gasch\"utz}\label{suga}

Next we need to report on a result of Gasch\"utz (see \cite{G}, \cite{Gr} and
\cite{JR}). 
For this, let $G$ be a finite group and $\pi : F_n\to G$ a 
surjective homomorphism
of a rank $n$ free group $F_n$ onto $G$. Let $R$ be the kernel of 
$\pi$ and $\bar R:=R/R'$ the corresponding relation module with its action  
(from the left) of $G$. The result of Gasch\"utz says that there is 
a $\BQ[G]$-module isomorphism
\begin{equation}\label{gasch}
\kappa:\BQ\oplus \BQ[G]^{n-1}\to  \BQ\otimes_{\bb Z}\bar R.
\end{equation}
where $\BQ$ is to be the trivial one dimensional $\BQ[G]$-module. 
We fix a $\BQ[G]$-module isomorphism $\kappa$ once for all and identify
$\BQ\otimes_\BZ\bar R$ with $\BQ\oplus \BQ[G]^{n-1}$. In
particular we obtain  ${\rm dim}_\BQ(\BQ\otimes_\BZ\bar{R})=1+|G|(n-1)$,
which is also clear from the formula of Reidemeister and Schreier.
We deduce 
\begin{equation}\label{gaschi1}
{\bf I}_1(\BQ\otimes_{\bb Z}\bar R)={\bf I}_1(\BQ[G])^{n}=\BQ^n,
\end{equation}
\begin{equation}\label{gaschi}
{\bf I}_i(\BQ\otimes_{\bb Z}\bar R)={\bf I}_i(\BQ[G])^{n-1}
\qquad (i=2,\ldots,\ell).
\end{equation}

From (\ref{dim17}) we get
\begin{lemma}\label{dim19}
Let $G$ be a finite group and $\pi : F_n\to G$ a 
surjective homomorphism. Let $\bar R:=R/R'$ be the corresponding 
relation module. Then we have:
$${\rm dim}_{\bb Q}({\rm Hom}_G({\bf I}_1(\BQ\otimes_{\bb Z}\bar{R}), 
{\bf I}_1(\BQ[G])))=n$$
$${\rm dim}_{\bb Q}(
{\rm Hom}_G({\bf I}_i(\BQ\otimes_{\bb Z}\bar{R}), {\bf I}_i(\BQ[G])))=
(n-1){\rm dim}_{\bb Q}({\bf I}_i(\BQ[G])$$
\end{lemma}

The result of Gasch\"utz enables us further to decompose the algebraic group
${\cal G}_{G,\pi}$ and to describe the components. 
We define for $i=1,\ldots,\ell$:
\begin{equation}\label{stabi}
{\cal G}_{G,\pi,i}:=\left\{\, \phi\in {\cal G}_{G,\pi}
\Mid \phi |_{{\bb C}\otimes{\bf I}_j({\bb Q}\otimes_{\bb Z}\bar R)}=
{\rm Id}_{{\bb C}\otimes{\bf I}_j({\bb Q}\otimes_{\bb Z}\bar R)}\ 
{\rm for}\ j\ne i\, \right\}.
\end{equation}
Note that the groups ${\cal G}_{G,\pi,i}$ centralise each other and
have pairwise trivial intersection. It is also easy to see that we have an
internal product decomposition:
$$
{\cal G}_{G,\pi}=\prod_{i=1}^\ell\, {\cal G}_{G,\pi,i}
$$
of $\BQ$-defined algebraic groups.

We write $\GL(d,S)$ ($d\in \BN$) for the group of 
invertible $d\times d$-matrices with entries in an associative ring $S$.
Following the discussion of the previous section we infer that the
identification made using  
(\ref{gasch}) defines isomorphisms
$$
\Sigma_i : {\cal G}_{G,\pi,i}(\BQ)\to \GL(n-1,{\bf I}_i(\BQ[G])^{\rm op})
=\GL((n-1)h_i,D_i^{\rm op})
$$
for $i=2,\dots \ell$ and also an isomorphism
$$
\Sigma_1 : {\cal G}_{G,\pi,1}(\BQ)\to \GL(n,{\bf I}_1(\BQ[G])^{\rm op})=
\GL(n,\BQ).
$$
Together they define an isomorphism
$$
\Sigma : {\cal G}_{G,\pi}(\BQ)\to \GL(n,\BQ)\times \prod_{i=2}^\ell 
\GL((n-1)h_i,D_i^{\rm op}).
$$

Suppose that $i$ is one of the numbers $1,\ldots,\ell$ and let 
$K_i$ be the center of ${\bf I}_i(\BQ[G])^{\rm op}$. As already
mentioned, this is a cyclotomic number field. The group 
$\GL((n-1)h_i,D_i^{\rm op})$ is the group of $K_i$-rational points of a
$K_i$-defined linear algebraic group ${\cal H}_i$. The above discussion shows 
that the $\BQ$-defined algebraic group ${\cal G}_{G,\pi,i}$ 
is equal to the base
field restriction of ${\cal H}_i$ from $K_i$ to $\BQ$. For this concept see
\cite{PLR}. We shall add a description of the 
groups ${\cal G}_{G,\pi}^1$ introduced in the introduction.
Every element $g$ of $\GL((n-1)h_i,D_i^{\rm op})$ ($i=1,\ldots,\ell$) 
acts by multiplication on $D_i^{(n-1)h_i}$. Taking a $K_i$-vector space
identification $D_i=K_i^{d_i}$ we associate to $g$ a square matrix of size 
$(n-1)h_id_i$. The determinant of this matrix is called the reduced norm 
${\rm Nrd}(g)$ of $g$. The group of elements of reduced norm $1$ is
conventionally denoted by $\SL((n-1)h_i,D_i^{\rm op})$.
We have
$$
{\cal G}^1_{G,\pi,i}(\BQ)=\{\, g\in {\cal G}_{G,\pi,i}(\BQ)\Mid 
{\rm Nrd}(\Sigma_i(g))=1\,\}.
$$

\subsection{Representations of subgroups of finite index in
  $\A(F_n)$}\label{repsfini}

In this section we describe a matrix version of the representation (\ref{I5}).
To do this we use Lemma \ref{repkon} which constructs an integral 
matrix representation of a subgroup of finite index in $\Gamma(G,\pi)$ which
is compatible with the representations
$\Sigma_i$ ($i=1,\ldots,\ell$) and $\Sigma$
from the previous section. 

Let $G$ be a finite group and $\pi : F_n\to G$ a surjective homomorphism
of the free group $F_n$ onto $G$. Let $R$ be the kernel of 
$\pi$ and $\bar R$ the corresponding relation module. Let further
$\Gamma(G,\pi)$ be the subgroup of $\A(F_n)$ defined in 
(\ref{I1}). In the sequel we shall often use the following fact without
further notice.
\begin{lemma} An element $\varphi\in\Gamma(R)$ is in $\Gamma(G,\pi)$ if and
only if the map $\bar\varphi :\bar R\to \bar R$ induced by $\varphi$ is
$G$-equivariant. 
\end{lemma}
\begin{proof} If $\varphi$ is in $\Gamma(G,\pi)$ then  
for every $w\in F_n$ we have a $r_w\in R$ with $\varphi(w)=wr_w$. Hence
$\varphi(wrw^{-1})=wr_w\varphi(r)r_w^{-1}w^{-1}$ holds for every $r\in R$.
The latter is equal to $ w\varphi(r)w^{-1}$ modulo $R'$ which 
implies that $\bar\varphi :\bar R\to \bar R$ is
$G$-equivariant.

If $\bar\varphi :\bar R\to \bar R$ is $G$-equivariant then for every 
$w\in F_n$ the congruence
$\varphi(w)\varphi(r)\varphi(w)^{-1}\equiv w \varphi(r) w^{-1}$ 
holds modulo
$R'$ for every $r\in R$. This means that $w$ and $\varphi(w)$ act the same on
$\bar R$. The result of Gasch\"utz implies that the action of $G$ on $\bar R$
is faithful, that is no element except the identity of $G$ acts as the
identity. We conclude that $\varphi(w)w^{-1}$ is in $R$.
\end{proof}

Keeping the notations set up in subsections \ref{sura}, \ref{suga} we define  
$$
\Lambda_i:=\bar R\cap {\bf I}_i(\BQ\otimes_{\bb Z} \bar R)
\quad (i=1,\ldots ,\ell).
$$
Each $\Lambda_i$ is a $\BZ[G]$-submodule of 
$\bar R\cap {\bf I}_i(\BQ\otimes_{\bb Z} \bar R)$ and 
$\prod_{i=1}^\ell \Lambda_i\le \bar R$
is a finite index $\BZ[G]$-submodule.
Notice that the $\Lambda_i$ is a sublattice of 
${\bf I}_i(\BQ\otimes_{\bb Z} \bar R)=M(h_i,D_i)^{n-1}$ 
($i=2,\ldots,\ell$) in the
sense of the previous section, so is 
$\Lambda_1$ in ${\bf I}_1(\BQ\otimes_{\bb Z} \bar R)=\BQ^{n}$. Applying Lemma
\ref{repkon} we choose orders ${\cal R}_i\subset D_i$ and finite index
additive subgroups $\tilde\Lambda_i\le \Lambda_i$ ($i=1,\ldots,\ell$). 
We define 
\begin{equation}\label{hut}
\hat\Gamma(G,\pi)\le \Gamma(G,\pi)
\end{equation}
to be the simultaneous stabiliser of all the 
$\tilde\Lambda_i\le \Lambda_i$ ($i=1,\ldots,\ell$).
Notice further that $\hat\Gamma(G,\pi)$
is of finite index in  $\Gamma(G,\pi)$. Using the various representations
$\Theta$ from Lemma \ref{repkon} we obtain representations 
$$
\sigma_{G,\pi,1}: \hat\Gamma (G,\pi) \to \GL(n,\BZ),\qquad
\sigma_{G,\pi,i}: \hat\Gamma (G,\pi) \to 
\GL((n-1)h_i,{\cal R}_i^{\rm op})
$$
($i=2,\ldots, \ell$), 
which come by projection from the analogously defined representation 
$$
\sigma_{G,\pi}: \hat\Gamma(G,\pi) \to 
\GL(n,\BZ)\times \prod_{i=2}^\ell \GL(n-1,M(h_i,{\cal R}_i^{\rm op})).
$$
Notice that the diagram 
$$
\xymatrix{
\hat\Gamma(G,\pi) \ar[r]^{} \ar[d]^{\sigma_{G,\pi}}
& {\cal G}_{G,\pi}(\BQ) \ar[ld]^{\Sigma} \\
\GL(n,\BQ)\times \prod_{i=2}^\ell \GL(n-1,M(h_i,D_i^{\rm op})) & 
}
$$
is commutative.

\section{Derivations and the relation module}\label{Drei}

This section contains a brief discussion of the 
correspondences between Fox-derivatives, the
relation module of a finite group and cohomology. We apply these results to
obtain a certain canonical system of right module generators of 
${\rm Hom}_G(\BQ\otimes_{\bb Z} \bar R,\BQ[G])$ 
where $G$ is a finite group and 
$\bar R$ is the relation module coming from a presentaion $\pi :F_n\to G$.

\subsection{Fox calculus and cohomology}

For $n\in\BN$, let $F_n$ be the free group with basis $x_1,\ldots ,x_n$. We
write $\BZ[F_n]$, $\BQ[F_n]$ for the group rings of $F_n$ over the integers or
the rational numbers, respectively.
Let 
$$
\frac{\partial}{\partial x_i}:F_n\to \BZ[F_n]\qquad (i=1,\ldots n)
$$
be the Fox-derivatives (see \cite{Bi} for a more elaborate exposition). 
They are defined by the rules
\begin{equation}\label{F1}
\frac{\partial x_j}{\partial x_i}=\delta_{ij},\qquad \ 
\frac{\partial w_1w_2}{\partial x_i}=\frac{\partial w_1}{\partial x_i}
+w_1\frac{\partial w_2}{\partial x_i}
\end{equation}
for $i,\, j=1,\ldots n$ and $\ w_1,\, w_2\in F_n$. As usual,
$\delta_{ij}$ denoting the Kronecker symbol ($\delta_{ij}=1$ if $i=j$ and 
$\delta_{ij}=0$ if $i\ne j$). For later use we note some consequences of 
(\ref{F1}). For $i=1,\ldots , n$ and $e\in \BN$ we have 
\begin{equation}\label{F2}
\frac{\partial x_i^e}{\partial x_i}=\sum_{k=0}^{e-1}x_i^k,\qquad\quad 
\frac{\partial x_i^{-e}}{\partial x_i}=-\sum_{k=1}^{e}x_i^{-k}.
\end{equation}
Let further $w:=v_1x_i^{e_1} v_2x_i^{e_2}\cdot 
\ldots \cdot v_kx_i^{e_k}v_{k+1}$ ($e_1,\ldots ,e_k\in\BZ$)
be an element of $F_n$ with the subwords 
$v_j$ ($j=1,\ldots ,k+1$) not involving $x_i$, then
\begin{equation}\label{F3}
\frac{\partial w}{\partial x_i}=v_1\frac{\partial x_i^{e_1}}{\partial x_i}+
v_1x_i^{e_1}v_2 \frac{\partial x_i^{e_2}}{\partial x_i}+\ldots +
v_1x_i^{e_1}\cdot v_2x_i^{e_2}\cdot \ldots \cdot v_k
\frac{\partial x_i^{e_k}}{\partial x_i}.
\end{equation}
Let $G$ be a finite group generated by its elements
$g_1,\ldots , g_n$ ($n\in \BN$) and let
$$
\pi: F_n\to G\qquad \pi(x_1)=g_1,\ldots ,\pi(x_n)=g_n
$$
be the corresponding surjective group homomorphism.  
The homomorphism $\pi:F_n\to G$ extends by linearity to 
surjective ring homomorphisms:
$$
\pi: \BZ [F_n]\to \BZ[G]\qquad \pi: \BQ [F_n]\to \BQ[G].
$$
Let $R\le F_n$ be the kernel of $\pi$ and set $\bar{R}:=R/R'$
to be the commutator factor group of $R$. Let $\bar{} :R\to R/R'$ stand for
the quotient homomorphism.
The free group $F_n$ acts by
conjugation on $R$. Since its subgroup $R$ acts trivially on 
$\bar{R}$ we get an induced action of $G$ on $\bar{R}$ which
we denote by $g\cdot u$ ($g\in G$, $u\in \bar{R}$). We generally write the
composition in $\bar{R}$ additively.   
The free abelian group $\bar{R}$ then obtains the structure of a left
$\BZ[G]$-module which extends to a left $\BQ[G]$-module structure on
$\BQ\otimes_\BZ\bar{R}$. We have
\begin{lemma}\label{leF0}
Let $G$ be a group and $\pi :F_n\to G$ a homomorphism. 
Let $R$ be the kernel of $\pi$. Then
$$\pi\left( \frac{\partial r}{\partial x_i}\right)=0$$
for all $r\in R'$ and $i=1,\ldots n$.
\end{lemma}
\begin{proof}
Note first that (\ref{F1}) implies 
$\frac{\partial w^{-1}}{\partial x_i}=-w^{-1}
\frac{\partial w}{\partial x_i}$ for all $w\in F_n$ and $i=1,\ldots n$.
It is enough to show the result for single commutators 
$r=r_1r_2r_1^{-1}r_2^{-1}$. A repeated application of (\ref{F1}) implies that 
$$\frac{\partial r}{\partial x_i}= 
\frac{\partial r_1}{\partial x_i}+
r_1\frac{\partial r_2}{\partial x_i}-
r_1r_2r_1^{-1}\frac{\partial r_1}{\partial x_i}-
r_1r_2r_1^{-1}r_2^{-1}\frac{\partial r_2}{\partial x_i}$$
for all $i=1,\ldots n$. As $\pi(r_1)=\pi(r_2)=1$ the result follows.
\end{proof}

The Fox-derivatives  
$\frac{\partial}{\partial x_i}$ now induce maps
$$
\partial_i: \bar{R}\to \BZ[G]\qquad (i=1,\ldots ,n)
$$
defined as follows.
For every $u\in\bar{R}$ choose $w_u\in R$ with $ \bar w_u=u$ and
define 
\begin{equation}\label{deldef} 
\partial_i(u):=\pi\left( \frac{\partial w_u}{\partial x_i}\right)\qquad
(i=1,\ldots ,n). 
\end{equation}
We have
\begin{lemma}\label{leF}
The maps $\partial_1,\ldots \partial_n$ do not depend on the choices of the 
preimages $w_u$ ($u\in\bar{R}$). They are $\BZ[G]$-module
homomorphisms when $\BZ[G]$ acts on itself from the left.
\end{lemma}
\begin{proof}
The first statement follows from Lemma \ref{leF0}. For the second note that
$$\frac{\partial grg^{-1}}{\partial x_i}=
\frac{\partial g}{\partial x_i}+g \frac{\partial r}{\partial x_i}-
grg^{-1}\frac{\partial g}{\partial x_i}$$
for all $g,\, r\in F_n$ and all $i=1,\ldots n$. 
When $r$ is in $R$, $\pi(r)=1$ and the projection of the right hand side into
$\BZ[G]$ is equal to $\pi(g) \pi(\frac{\partial r}{\partial x_i})$.
\end{proof}
We also write
$\partial_i: \BQ\otimes_{\bb Z}\bar{R}\to \BQ[G]$
for the induced $\BQ[G]$-module homomorphisms.

Let $t$ be the $\BZ$-rank of the free abelian group $\bar{R}$.
Note that $t=1+|G|(n-1)$ and also
$t={\rm dim}_{\bb Q}(\BQ\otimes_{\bb Z}\bar{R})$. 
Let
$B:=(u_1,\ldots, u_t)$ be a $\BQ$-basis of $\BQ\otimes_{\bb Z}\bar{R}$. 
We consider here
\begin{equation}\label{pm}
J_B:= \left( \begin{array}{ccc}
                \partial_1(u_1) & \ldots & \partial_n(u_1)\\
                \vdots &  & \vdots \\
                \partial_1(u_t) & \ldots & \partial_n(u_t)
\end{array} \right),\qquad
\end{equation}
which is a $t\times n$ matrix with entries in $\BQ[G]$.

Given a left $\BQ[G]$-module $M$ and $s\in\BN$ we write
$$M^s:=\left\{ 
\left(
\begin{array}{c}
m_1 \\ \vdots \\ m_s
\end{array} \right) \Mid m_1,\ldots, m_s\in M\right\}
$$
for the corresponding $\BQ[G]$-module of column vectors.
The matrix $J_B$ induces a $\BQ$-linear map
$J_B:M^n\to M^t$. Notice for later use that given a $\BQ[G]$-submodule
$N\le M$ we have $J_B(N^n)\subseteq N^t$.

Given the left $\BQ[G]$-module $M$ we write
$$
{\rm Der}(G,M):=\{\, d:G\to M\Mid d(gh)=gd(h)+d(g) \ {\rm for\ all}\ g,h\in
G\,\} 
$$
for the vector space of derivations from $G$ to $M$.
Using the generators $g_1,\ldots ,g_n$ of our finite group $G$ we define 
for $d\in {\rm Der}(G,M)$ 
$$
L(d):=
\left(\begin{array}{c}
d(g_1) \\ \vdots \\ d(g_n)
\end{array}\right)\in M^n.
$$
We have
\begin{lemma}\label{coholemm}
The map $L: {\rm Der}(G,M)\to M^n$ is injective and its image is equal to
${\rm Ker}(J_B)$. Hence $L$ defines a $\BQ$-vector space isomorphism
between ${\rm Der}(G,M)$ and ${\rm Ker}(J_B)$.
\end{lemma}
\begin{proof} 
Let $d$ be a derivation in the kernel of $L$. We infer that
$d(g_1)=\ldots=d(g_n)=0$. The derivation rule and the fact that the 
$g_1,\ldots ,g_n$ generate $G$ implies that $d$ is identically zero.
This shows that $L$ is injective.

Let $D:F_n\to M$ be a derivation where $F_n$ acts on $M$ through the
homomorphism $\pi : F_n\to G$. Since both sides are derivations and agree on
the $x_i$'s we find 
\begin{equation}\label{deriv}
D(w)=\frac{\partial w}{\partial x_1}D(x_1)+\ldots
+\frac{\partial w}{\partial x_n}D(x_n)
\end{equation}
for all $w\in F_n$. 

Given a derivation $d: G\to M$, define an induced derivation 
$D: F_n\to M$ by setting $D(w):=d(\pi(w))$. Applying (\ref{deriv}) we find
for every $r\in R$
$$0=d(1)=D(r)=\sum_{i=1}^n \frac{\partial r}{\partial x_i}D(x_i)=
\sum_{i=1}^n \partial_i (r) d(g_i).$$
Applying this to $r=u_1,\ldots, u_t$ shows that the image of $L$ is contained
in the kernel of $J_B$.

Let now 
$$m:=\left(\begin{array}{c}
m_1 \\ \vdots \\ m_n
\end{array}\right)\in M^n
$$
be an $n$-tuple of elements of $M$. It is well known that there is a
derivation $D :F_n\to M$ with $D(x_1)=m_1,\ldots, D(x_n)=m_n$. 
If $m$ is an element of the kernel of $J_B$ then (\ref{deriv}) implies that
$D(r)=0$ for all $r\in R$. The derivation $D$ then induces a derivation
$d: G\to M$ with the property $d(g_1)=m_1,\ldots, d(g_n)=m_n$. This shows that
every element of the kernel of $J_B$ is in the image of $L$.
\end{proof}
\begin{corollary}\label{cohocor} 
Let $M$ be a $\BQ[G]$-module then the kernel of the map
$J_B :M^n\to M^t$ is zero if $M$ is a trivial module and has $\BQ$-vector
space dimension equal to ${\rm dim}_{\bb Q}(M)$ if M does not have a trivial
submodule. 
\end{corollary}
\begin{proof} 
Given an element $m\in M$ we let $d_m :G\to M$ be the corresponding inner
derivation $d_m(g)=(g-1)m$. Let  
$$
{\rm IDer}(G,M):=\{\, d_m \Mid m\in M\,\}\subseteq  {\rm Der}(G,M)
$$
be the space of inner derivations. It is well known that
the first cohomology group $H^1(G,M)={\rm Der}(G,M)/{\rm IDer}(G,M)$
is zero ($G$ is a finite group). Note also that the map 
$\mu : M\to {\rm IDer}(G,M)$, 
defined by $\mu(m):=d_m$ ($m\in M$) is surjective.

If $M$ is a trivial module then $\mu$ is identically zero, that implies  
${\rm IDer}(G,M)=\{0\}$ and hence also ${\rm Der}(G,M)=\{0\}$. 
We find from Lemma \ref{coholemm} that 
${\rm dim}_{\bb Q}({\rm ker}(J_B))=0.$

If $M$ does not contain the trivial module then 
$\mu$ is injective, that implies
${\rm dim}_{\bb Q}({\rm IDer}(G,M))={\rm dim}_{\bb Q}(M)$. We find from 
Lemma \ref{coholemm} that 
${\rm dim}_{\bb Q}({\rm ker}(J_B))={\rm dim}_{\bb Q}({\rm Der}(G,M))={\rm
  dim}_{\bb Q}(M).$
\end{proof}

The results of this section will be applied in the next section in case of the
relation modules.

\subsection{${\rm Hom}_G(\BQ\otimes_{\bb Z} \bar R,\BQ[G])$ 
as right $\BQ[G]$-module}\label{Hom}

We keep here our notation from before, that is $G$ is a finite group 
generated by its elements $g_1,\ldots , g_n$ ($n\in \BN$) and 
$\pi: F_n\to G$ is the homomorphism defined by 
$\pi(x_1)=g_1,\ldots ,\pi(x_n)=g_n$. Furthermore $R$ is the kernel of $\pi$
and $\bar R$ is the corresponding relation module.
As the main result we shall show that the 
$\partial_1,\ldots,\partial_n\in {\rm Hom}_G(\BQ\otimes_\BZ \bar R,\BQ[G])$,
which were constructed in the previous subsection generate  
${\rm Hom}_G(\BQ\otimes_{\bb Z} \bar R,\BQ[G])$ if this 
space is considered as a
$\BQ[G]$-module from the right in the following way. 

Let $M$ be a left $\BQ[G]$-module. For 
$f\in {\rm Hom}_G(M ,\BQ[G])$
and $A\in \BQ[G]$ we define $f^A\in {\rm Hom}_G(M ,\BQ[G])$ by 
\begin{equation}
f^A(x):=f(x)A\qquad (x\in M).
\end{equation}
Thus the $\BQ$-vector space ${\rm Hom}_G(M ,\BQ[G])$
becomes a right $\BQ[G]$-module.

We shall prove here:
\begin{proposition}\label{wich} The maps $\partial_1,\ldots, \partial_n$ 
(defined in (\ref{deldef})) generate 
${\rm Hom}_G(\BQ\otimes_{\bb Z} \bar{R},\BQ[G])$ 
as right $\BQ[G]$-module. 
\end{proposition}

To prepare for the proof we decompose the matrix $J_B$ into submatrices. 
To do this we choose
the basis $B$ adapted to the decomposition
$$
\BQ\otimes_{\bb Z}\bar{R}=\prod_{i=1}^\ell\, 
{\bf I}_i(\BQ\otimes_{\bb Z}\bar{R}).
$$
That is we choose bases
$B_i=(u_{i1},\ldots u_{it_i})$ 
$(t_i:={\rm dim}_{\bb Q}({\bf I}_i(\BQ\otimes_{\bb Z}\bar{R}))
$
for $i=1,\ldots,\ell$ and put them together to obtain a 
$\BQ$-vector space basis 
$$
B=(u_1=u_{11},\ldots,u_t=u_{\ell t_\ell}) 
$$
of $\BQ\otimes_{\bb Z}\bar{R}$.

Consider the matrices 
\begin{equation}\label{pmi}
J_B^{P_i}:= \left( \begin{array}{ccc}
                \partial_1(u_1)P_i & \ldots & \partial_n(u_1)P_i\\
                \vdots &  & \vdots \\
                \partial_1(u_t)P_i & \ldots & \partial_n(u_t)P_i
\end{array} \right)=
\left( \begin{array}{ccc}
                \partial_1^{P_i}(u_1) & \ldots & \partial_n^{P_i}(u_1)\\
                \vdots &  & \vdots \\
                \partial_1^{P_i}(u_t) & \ldots & \partial_n^{P_i}(u_t)
\end{array} \right)
\end{equation}
for $i=1,\ldots,\ell$. Notice that according to 
property {\bf LP} of the right projector $P_i$ (see section \ref{sura})
the maps
$$J_B,\, J_B^{P_i}: {\bf I}_i(\BQ[G])^n\to {\bf I}_i(\BQ[G])^t$$
are the same. Furthermore we have 
$\partial_j(u)\in {\bf I}_i(\BQ[G])$ ($j=1,\ldots ,n$)
for $u\in {\bf I}_i(\BQ\otimes_{\bb Z}\bar{R})$. 
This follows since each $\partial_j$ is a $\BQ[G]$-module homomorphism which
has to respect isotypic components.
Thus property {\bf RP} of the right projector $P_i$
implies that all entries of $J_B^{P_i}$ are equal to zero except
those in the submatrix
\begin{equation}\label{pmis}
J_{B_i}:= \left(\begin{array}{ccc}
                \partial_1(u_{i1})P_i & \ldots & \partial_n(u_{i1})P_i\\
                \vdots &  & \vdots \\
                \partial_1(u_{i t_i})P_i & \ldots & \partial_n(u_{i t_i})P_i
\end{array} \right)
\end{equation}
hence the kernels of the two maps 
$J_B : {\bf I}_i(\BQ[G])^n\to {\bf I}_i(\BQ[G])^t$ and
$J_{B_i}:{\bf I}_i(\BQ[G])^n\to {\bf I}_i(\BQ[G])^{t_i}$
are the same. Note finally that we have 
$$
J_{B_i}=\left(\begin{array}{ccc}
           \partial_1^{P_i}(u_{i1}) & \ldots & \partial_n^{P_i}(u_{i1})\\
           \vdots &  & \vdots \\
           \partial_1^{P_i}(u_{i t_i}) & \ldots & \partial_n^{P_i}(u_{i t_i})
\end{array} \right)
$$
for all $i=1,\ldots,\ell$.

\medskip
\begin{prof} {\it of Proposition \ref{wich}:}
We have the decomposition
$$
{\rm Hom}_G(\BQ\otimes_{\bb Z}\bar{R}, \BQ[G])=
\prod_{i=1}^\ell
{\rm Hom}_G({\bf I}_i(\BQ\otimes_{\bb Z}\bar{R}), {\bf I}_i(\BQ[G]))
$$
of ${\rm Hom}_G(\BQ\otimes_{\bb Z}\bar{R}, \BQ[G])$ as 
a right $\BQ[G]$-module. 
Since the $\partial_1^{P_i},\ldots ,\partial_n^{P_i}$ are zero on  
${\bf I}_j(\BQ\otimes_{\bb Z}\bar{R})$ for $j\ne i$ it is enough to show that
they generate 
${\rm Hom}_G({\bf I}_i(\BQ\otimes_{\bb Z}\bar{R}), {\bf I}_i(\BQ[G]))$ 
as a right $\BQ[G]$-module.
Consider the $\BQ$-linear map
$\Phi_i : {\bf I}_i(\BQ[G])^n \to 
{\rm Hom}_G({\bf I}_i(\BQ\otimes_{\bb Z}\bar{R}), {\bf I}_i(\BQ[G]))$
given by 
$$
\Phi_i\left(\left(
\begin{array}{c}
m_1 \\ \vdots \\ m_n
\end{array} \right)\right):=
\partial_1^{P_i m_1}+ \ldots +\partial_n^{P_i m_n}.
$$
Notice that ${\bf I}_i(\BQ[G])\cdot \BQ[G]={\bf I}_i(\BQ[G])$ and hence
the image of $\Phi_i$ is the right $\BQ[G]$-submodule of 
${\rm Hom}_G({\bf I}_i(\BQ\otimes_{\bb Z}\bar{R}), {\bf I}_i(\BQ[G]))$ 
generated by
the $\partial_1^{P_i},\ldots ,\partial_n^{P_i}$.
The key observation now is 
\begin{equation}\label{dime}
{\rm kernel}(J_{B_i})={\rm kernel}(\Phi_i)
\end{equation}
which is obvious by what is explained above. Notice that in fact 
$(m_1,\ldots,m_n)^t$ is in the kernel of $\Phi_i$ if and only if
$\partial_1^{P_im_1}(u_{ik})+\ldots +\partial_n^{P_im_n}(u_{ik})=0$
holds for all $k=1,\ldots,t_i$.

We can treat now the case $i=1$ (the case of the trivial module).
Here we have by (\ref{gaschi1})  that  
${\rm dim}_{\bb Q}({\bf I}_1(\BQ\otimes_{\bb Z}\bar{R}))=n$ and the kernel 
of $\Phi_1$ is trivial (by Corollary \ref{cohocor}). Since we have
${\rm dim}_{\bb Q}({\rm Hom}_G({\bf I}_1(\BQ\otimes_{\bb Z}\bar{R}), 
{\bf I}_1(\BQ[G])))=n$ (by Lemma \ref{dim19}) the map $\Phi_1$ has to be
surjective. 

In the cases $i >1$ we have
$ {\rm dim}_{\bb Q}({\rm kernel}(\Phi_i))={\rm dim}_{\bb Q}({\bf I}_i(\BQ[G]))$
by (\ref{dime}) and Corollary \ref{cohocor}. By Lemma \ref{dim19} and 
(\ref{dim17}) we have
${\rm dim}_{\bb Q}(
{\rm Hom}_G({\bf I}_i(\BQ\otimes_{\bb Z}\bar{R}), {\bf I}_i(\BQ[G]))=
(n-1){\rm dim}_{\bb Q}({\bf I}_i(\BQ[G]))$
which implies that $\Phi_i$ has to be surjective.
\end{prof}

As a corollary to Proposition \ref{wich} we note
\begin{corollary}\label{wichcor} 
The maps 
$\partial_1^{P_i},\ldots ,\partial_n^{P_i}$  
generate 
${\rm Hom}_G({\bf I}_i(\BQ\otimes_{\bb Z}\bar{R}), {\bf I}_i(\BQ[G]))$ 
as right ${\bf I}_i(\BQ[G])$-module for every $i=1,\ldots ,\ell$. 
\end{corollary}
We shall give now a reformulation of the last corollary which will be needed
later. Let  
${\rm Id}_{i,k}\in {\bf I}_i(\BQ\otimes_{\bb Z}\bar{R})$ 
($k=1,\ldots,n-1$) be the images 
under $\kappa_i$ (see (\ref{gaschi}), (\ref{gaschi1})) 
of the unit elements of of the components of 
${\bf I}_i(\BQ\otimes_{\bb Z}\bar{R})^{n-1}$ respectively of
${\bf I}_1(\BQ\otimes_{\bb Z}\bar{R})^{n}$ with $k=1,\ldots,n$. 
Using the identification (\ref{ra1oppp}) we think of  
$f\in{\rm Hom}_G({\bf I}_i(\BQ\otimes_{\bb Z}\bar{R}), {\bf I}_i(\BQ[G]))$ 
being identified with the collection of matrices 
$f({\rm Id}_{i,1})=A_{i,1},\ldots ,
f({\rm Id}_{i,n-1})=A_{i,n-1}\in M(h_i,D_i^\op)$ or 
$f({\rm Id}_{1,1})=A_{1,1},\ldots ,
f({\rm Id}_{i,n})=A_{1,n}\in M(h_1,D_1^\op)=\BQ$.

\begin{corollary}\label{wichcor1}
Let $i$ be one of the numbers $2,\ldots,\ell$.
Given $A_{1},\ldots ,A_{n-1}\in M(h_i,D_i^\op)$
there are $B_1,\ldots,B_n\in \BQ[G]$ such that
$$\partial_1^{P_iB_1}+\ldots +\partial_n^{P_iB_n}({\rm Id}_{i,k})=A_k
\qquad (k=1,\ldots,n-1).$$
Given $A_{1},\ldots ,A_{n}\in M(h_1,D_1^\op)$
there are $B_1,\ldots,B_n\in \BQ[G]$ such that
$$\partial_1^{P_1B_1}+\ldots +\partial_n^{P_1B_n}({\rm Id}_{1,k})=A_k
\qquad (k=1,\ldots,n).$$
\end{corollary}

{\it Remark:} We sketch here another (more conceptual) 
proof of Proposition \ref{wich}. Features of the proof given above 
will play a role in the sequel. 
 
Using $H^2(G,\BQ\otimes_{\bb Z}\bar R)=0$ we 
obtain an injective homomorphism 
$$\epsilon : F_n/R'\to (\BQ\otimes_{\bb Z}\bar R)\rfish G$$
from $F_n/R'$ to the semidirect product on the right hand side. 
On $R/R'\le F_n/R'$ the homomorphism $\epsilon$ induces the inclusion of 
$R/R'$ into $\BQ\otimes_{\bb Z}\bar R\le (\BQ\otimes_{\bb Z}\bar R)\rfish G$.
The second component of $\epsilon$ coincides with $\pi$.  
Associate to $f\in{\rm Hom}_G(\BQ\otimes_{\bb Z}\bar R,\BQ[G])$ the 
homomorphism $\tau_f: (\BQ\otimes_{\bb Z}\bar R)\rfish G\to \BQ[G]\rfish G$
($\tau_f((u,g))=(f(u),g)$). Composition with $\epsilon$ creates a homomorphism
$$\lambda_f: F_n\to F_n/R'\to \BQ[G]\rfish G, \qquad
\lambda_f(w)=(d_f(w),\pi(w)),\quad (w\in F_n)$$
The map $d_f:F_n \to \BQ[G]$ is a derivation. Hence we have for $w\in F_n$
$$d_f(w)=\frac{\partial w}{\partial x_1}d_f(x_1)+\ldots
+\frac{\partial w}{\partial x_n}d_f(x_n).
$$
Restricting this formula to $R/R'$ we have expressed $f$ in the desired way.

\section{Redundant presentations}\label{Vier}

In this section, we present a special decomposition of the relation module
$\bar R$ for redundant presentations of finite groups. We then describe the
action of various Nielsen automorphisms of $F_n$ on $\bar R$ with respect to
this decomposition. 

Let $F_{n}$ ($n\in\BN, n\ge 2$)
be the free group with basis $x_1,\ldots ,x_{n-1},\, y$.
Let $G$ be a finite group generated by its elements
$g_1,\ldots , g_{n-1}$  and let
\begin{equation}\label{redu1}
\pi: F_{n}\to G\qquad \pi(x_1)=g_1,\ldots ,\pi(x_{n-1})=g_{n-1},\ \pi(y)=1
\end{equation}
be the surjective group homomorphism corresponding to these data. 
We keep the notation introduced above of objects related to the presentation
$\pi$ of $G$. That is $R$ is the kernel of $\pi$ and $\bar R=R/R'$ is the
relation module, etc.. 
We shall construct certain Nielsen-like elements in $\hat\Gamma(G,\pi)$ 
(see \ref{hut}) which
will allow us to find many unipotent elements in the image of the
representation  
\begin{equation}\label{rep12}
\sigma_{G,\pi}: \hat\Gamma(G,\pi) \to 
\GL(n,\BZ)\times \prod_{i=2}^\ell \GL(n-1,M(h_i,{\cal R}_i^\op))
\end{equation}
defined in Section \ref{repsfini}.

\subsection{Decomposition of the relation module}

With the notation fixed as indicated above we define
$$
R_{\rm old}:=\{\, x\in \langle x_1,\ldots ,x_{n-1}\rangle\Mid \pi(x)=1\,\}\ 
{\rm and}\ \Ro:=R_{\rm old}/R'_{\rm old}.
$$
The injection $R_{\rm old}\to R$ defines a $\BZ[G]$-module homomorphism
$\Ro\to \bar R$. Since $y\in R$, we can consider the $\BZ[G]$-submodule
of $\bar R$ generated by $\bar y\in \bar R$. 

\begin{lemma} The homomorphism $\Ro\to\bar R$ is injective and so is the
homomorphism
$\BZ[G]\to \bar R$ defined by $A\mapsto A\cdot \bar y$.
Identifying $\Ro$ with its image in $\bar R$ we obtain a direct
sum decomposition
$$
\bar R=\bar R_{\rm old}\oplus \BZ[G]\cdot \bar y.
$$
as a $\BZ[G]$-module.
\end{lemma}
\begin{proof}
Let us temporarily write $\hat R_{\rm old}$ for the image of $R_{\rm old}$ in
$\bar R$. We first show that $\bar R=\hat R_{\rm old}+ \BZ[G]\cdot \bar y$.
Given any $w\in F_n$, by looking at the image of $w$ modulo the normal closure
of $y$ in $F_n$, we see that $w$ can be written in a unique way as 
$w=w_1w_2$ with $w_1,\, w_2\in F_n$ such
that $w_1$ does not involve the letter $y$ whereas $w_2$ is a product of
conjugates of powers of $y$. If $w$ is an element of $R$ then $w_1$ has to be
in $R_{\rm old}$. This implies $\bar R=\hat R_{\rm old}+ \BZ[G]\cdot \bar y$.

We infer that $\BQ\otimes_{\bb Z}\bar R=\BQ\otimes_{\bb Z}\hat R_{\rm old}
+ \BQ[G]\cdot \bar y$. We now count dimensions using the result of 
Gasch\"utz. As $\bar R$ and $\bar R_{\rm old}$ are free over $\BZ$, 
the lemma follows.
\end{proof}

\subsection{Nielsen Automorphisms}

We study here the effect of certain Nielsen automorphisms on the relation
module $\bar R$.
We define elements $\alpha_i\in \A(F_n)$ ($i=1,\ldots ,n-1$) by
\begin{equation}\label{al}
\alpha_i(x_i):=yx_i.
\end{equation}
Our convention here is that values not given are identical to the argument,
that is we assume  
$\alpha_i(x_j):=x_j,\ {\rm for}\ j\ne i\ {\rm and}\ \alpha_i(y)=y.$

For $X\in \langle x_1,\ldots ,x_{n-1}\rangle$ we define elements 
$\beta_X\in \A(F_n)$ by
\begin{equation}\label{be}
\beta_X(y):=XyX^{-1}
\end{equation}
and finally for $U\in  R_{\rm old}$ we define $\gamma_U\in \A(F_n)$ by
\begin{equation}\label{ge}
\gamma_U(y):=Uy.
\end{equation}
The $\alpha_i$, $\beta_X$, $\gamma_U$ extend to automorphisms of $F_{n}$.
We have
\begin{lemma}\label{NLe} 
The $\alpha_i$ ($i=1,\ldots,n-1$), $\beta_X$
($X\in \langle x_1,\ldots ,x_{n-1}\rangle$) and $\gamma_U$ 
($U\in  R_{\rm old}$) are in 
$\Gamma(G,\pi)$ and the following formulas hold for the 
corresponding automorphisms induced on $\bar R$:
\begin{equation}\label{alf}
\bar\alpha_i(u)=u+\partial_i(u)\cdot \bar y\quad {\it for}\ 
u\in \bar R_{\rm old}\ {\it and}\qquad
\bar\alpha_i(\bar y)=\bar y.
\end{equation}
\begin{equation}\label{bef}
\bar\beta_X(u)=u,\quad {\it for}\ u\in \bar R_{\rm old}\ {\it and}\qquad 
\bar\beta_X(\bar y)=\pi(X)\cdot \bar y,
\end{equation}
\begin{equation}\label{gef}
\bar\gamma_U(u)=u,\quad {\it for}\ u\in \bar R_{\rm old}\ {\it and}\qquad 
\bar\gamma_U(\bar y)=\bar U+\bar y,
\end{equation}
\end{lemma}
\begin{proof}
It is first of all clear from our definitions that all the $\alpha_i$, 
$\beta_X$, $\gamma_U$ are elements of\ 
$\Gamma(R):=\{\, \varphi\in \A(F_n)\Mid \varphi(R)=R \, \}$. Moreover all of
them induce the identity on $F_n/R$ and hence are even in $\Gamma(G,\pi)$. It
follows (see Section \ref{repsfini}) that they act 
as $\BZ[G]$-automorphisms on $\bar R$.

Next we prove formula (\ref{alf}). Fix $i\in\{1,\ldots ,n-1\}$ 
and define for $e\in\BN$
$$
W_e:=y\cdot x_iyx_i^{-1}\cdot\ldots \cdot x_i^{e-1}yx_i^{1-e},\quad
W_{-e}:=x_i^{-1}y^{-1}x_i\cdot\ldots \cdot x_i^{-e}y^{-1}x_i^{e}\in F_n 
$$
and add $W_0:=1$. We then have 
$\alpha_i(x_i^e)=W_e\cdot x_i^e$
for all $e\in \BZ$. For $u\in\Ro$ choose $w\in R_{\rm old}$ with $\bar w=u$.
We decompose $w\in \langle x_1,\ldots ,x_n\rangle$ as
$w:=v_1x_i^{e_1} v_2x_i^{e_2}\cdot 
\ldots \cdot v_kx_i^{e_k}v_{k+1}$
($e_1,\ldots ,e_k\in\BZ$) such that the $v_j$ do not involve $x_i$. Define
further 
$$
V_1:=v_1,\ V_2:=v_1x_i^{e_1}v_2,\ldots ,V_k:=v_1x_i^{e_1}v_2\ldots 
x_i^{e_k}v_k.
$$
We then have
$$
\alpha_i(w)=V_1W_{e_1}V_1^{-1}\cdot V_2W_{e_2}V_2^{-1}
\cdot\ldots\cdot 
V_k W_{e_2} V_k^{-1}\cdot w.
$$
Switching to the additive notation and using formulas (\ref{F2}) we find
$$
\bar\alpha_i(u)=u+\pi\left(v_1\frac{\partial x_i^{e_1}}{\partial x_i}+
v_1x_i^{e_1}v_2 \frac{\partial x_i^{e_2}}{\partial x_i}+\ldots +
v_1x_i^{e_1} v_2x_i^{e_2}\cdot \ldots \cdot v_k
\frac{\partial x_i^{e_k}}{\partial x_i}\right)\cdot \bar y
$$
and formula (\ref{alf}) follows from (\ref{F3}) and (\ref{deldef}). 

To prove formula (\ref{bef}) notice that $\beta_X$ is the identity on $\Ro$, 
subsequently (\ref{bef}) is clear.  
  
All statements made about the automorphisms $\gamma_U$ ($U\in R_{\rm old}$)
are evident. Notice that the map $\bar\gamma_U$ depends only on the image 
$\bar U$ of $U$ in $\bar R_{old}$.
\end{proof}
For later use we introduce now a notation for some elements in 
$\Gamma(G,\pi)$ obtained by composing the $\alpha_i$ and $\beta_X$. For
every $g\in G$ we choose $X_g\in \langle x_1,\ldots ,x_{n-1}\rangle$ 
with $\pi(X_g)=g$ and
define $\beta_g:=\beta_{X_g}$. We infer from formulas (\ref{alf}) and 
(\ref{bef})
$$\bar\beta_g\bar\alpha_i^m\bar\beta_{g^{-1}}(u+A\cdot\bar y)=
u+(\partial_i(u)mg+A)\cdot\bar y\qquad (u\in \Ro,\, A\in\BZ[G])
$$
for all $g\in G$, $m\in \BZ$.
Given an element $B=m_1g_1+\ldots +m_kg_k\in \BZ[G]$ 
($m_1,\ldots ,m_k\in\BZ$, $g_1,\ldots ,g_k\in G$) we define
$\eta_{B,i}\in\A(F_n)$ by
\begin{equation}\label{eta}
\eta_{B,i}:=\beta_{g_1}\alpha_i^{m_1}\beta_{g_1^{-1}}\circ\ldots \circ
\beta_{g_k}\alpha_i^{m_k}\beta_{g_k^{-1}}.
\end{equation}
for $i=1,\ldots ,n-1$.
By Lemma \ref{NLe} it follows that $\eta_{B,i}\in \Gamma(G,\pi)$ for
$i=1,\ldots ,n$ and $B\in \BZ[G]$.
Using the notation of Section \ref{Hom} we have
$$
\bar\eta_{B,i}(u+A\cdot\bar y)=u+(\partial_i^B(u)+A)\cdot\bar y
\qquad (u\in \Ro,\, A\in\BZ[G])
$$
for $i=1,\ldots ,n-1$ and $B\in\BZ[G]$. This formula makes it clear that
$$
\bar\eta_{m_1B_1+m_2B_2,i}=\bar\eta_{B_1,i}^{m_1}\circ \bar\eta_{B_2,i}^{m_2}
$$
holds for all $i=1,\ldots ,n-1$ and $B_1, \, B_2\in\BZ[G]$ and $m_1,\, m_2\in
\BZ$. 

Using the representation (\ref{furep}) we set
$$
\A^+(F_n):=\{\, \phi\in \A(F)\Mid {\rm det}(\rho_1(\phi))=1\, \}.
$$
The subgroup $\A^+(F_n)$ is of index $2$ in $\A(F_n)$. For $n\ge 3$ it is 
a perfect group, that is $\A^+(F_n)$ is equal to its own commutator subgroup. 
We note that all elements $\bar\eta_{B,i}$ are in $\A^+(F_n)$.

\section{Arithmeticity in case of a redundant presentation}\label{Fuenf}

In this section we prove Theorem \ref{teo}: The decomposition of $\bar R$ to
its isotypic components shows that our
main goal is to prove that the image of $\hat\Gamma(G,\pi)$ in (\ref{rep12})
contains a finite index subgroup of each component 
$$\SL(n-1,M(h_i,{\cal R}_i^\op))=\SL((n-1)h_i,{\cal R}_i^\op).$$ 
To prove this we will show that this image contains sufficiently
many unipotent elements with nonzero entries in the bottom strip of $h_i$ rows
and also in the most right hand strip of $h_i$ columns of this  
$\SL((n-1)h_i,{\cal R}_i^\op)$. Then we can appeal to a theorem of Vaserstein
which ensures that the elementary matrices with entries from a non-zero
two-sided ideal of ${\cal R}_i$ generate a finite index subgroup of the
arithmetic group $\SL((n-1)h_i,{\cal R}_i^\op)$.

\subsection{Elementary matrices}\label{Fuenf1}

This section sets up certain notations and facts about elementary matrices.

Let $S$ be an associative ring with unity.  
We denote by $\GL(n,S)$ 
($n\in\BN$) the group
of invertible $n\times n$ matrices with entries in $S$. We write $I_n$ for the
identity matrix. Given $1\le i,\, j\le n$ with $i\ne j$ and $s\in S$ we 
define $E_{ij}(s)\in \GL(n,S)$ to be the corresponding elementary matrix, that
is $E_{ij}(s)$ is equal to $I_n$ except for the $(i,\, j)$ position, where we
put $s$. For $1\le i,\, j\le n$ with $i\ne j$, $1\le k,\, l \le n$ 
with $k\ne l$ and $s_1,\, s_2\in S$
the Steinberg relations between the elementary matrices are 
\begin{equation}\label{stein}
[E_{ij}(s_1),E_{k\ell}(s_2)]:=\left\{\begin{array}{ccc}
                 1  &  {\rm for}\ j\ne k, & i\ne l, \\
                 E_{il}(s_1s_2) &  {\rm for}\ j= k, & i\ne l, \\
                 E_{kj}(-s_2s_1) &  {\rm for}\ j\ne k, & i= l, \\
\end{array} \right.
\end{equation}
where 
$[E_{ij}(s_1),E_{k l}(s_2)]=E_{ij}(s_1)E_{kl}(s_2)
E_{ij}(s_1)^{-1}E_{kl}(s_2)^{-1}$. For every admissible pair $i,\, j$ we have
$E_{ij}(s_1)E_{ij}(s_2)=E_{ij}(s_1+s_2)$ for all $s_1,\, s_2\in S$.

Let ${\lie a}\subset S$ be an ideal of $S$. Given natural numbers $m_1,\,
m_2\in\BN$ we denote by $M(m_1,m_2;{\lie a})$ the space of 
$m_1\times m_2$-matrices with entries in ${\lie a}$. 
Suppose now that $m_1+m_2=n$. For 
$A\in M(m_1,m_2;{\lie a})$ we define
$$
H(A):=\left(\begin{array}{cc}
           I_{m_2} &  0\\
           A &   I_{m_1}
\end{array} \right),\qquad
V(A):=\left(\begin{array}{cc}
           I_{m_2} &  A^t\\
           0 &   I_{m_1}
\end{array} \right)
$$
Here $A^t$ stands for the transpose of the matrix $A$. Suppose that 
$A:=(a_{ij})$ where $i=1,\ldots,m_1$, $j=1,\ldots,m_2$ with 
$a_{ij}\in{\lie a}$, then we have 
$$
H(A)=\prod_{i=1}^{m_1}   \prod_{j=1}^{m_2}\,     E_{m_2+i\, j}(a_{ij}),
\qquad 
V(A)=\prod_{i=1}^{m_1}   \prod_{j=1}^{m_2}\,     E_{j\, m_2+i}(a_{ij}).
$$
Given $m_1,\, m_2\in\BN$ with $m_1+m_2=n$ and an 
ideal $\lie{a}\le S$ we define 
$$
H(m_1,m_2;{\lie a}):=\{\, H(A)\Mid A\in M(m_1,m_2; {\lie a})\,\},
$$
$$
V(m_1,m_2;{\lie a}):=\{\, V(A)\Mid A\in M(m_1,m_2; {\lie a})\,\}.
$$
Both $H(m_1,m_2;{\lie a})$ and $V(m_1,m_2;{\lie a})$ 
are abelian subgroups of $\GL(n,S)$.

We apply this notation in case of the orders $M(h_i,{\cal R}_i)$ 
of the simple
factor rings ${\bf I}_i(\BQ[G])$ of the rational group ring $\BQ[G]$. 
So let $K$ be a number field and $D$ a 
finite dimensional division algebra over $K$ with center $K$. 
Let ${\cal R}$ be a order in $D$. It is well
known that every two-sided ideal $\lie{a}\le {\cal R}$ 
is additively finitely
generated and has finite index in ${\cal R}$ 
if it is non-zero. We have

\begin{proposition}\label{wasser} 
Let $n\in \BN$ with $n\ge 3$ and let $m_1,\, m_2\in\BN$ satisfy $m_1+m_2=n$.
Let ${\lie a}$ be a non-zero two-sided ideal of ${\cal R}$
Then $H(m_1,m_2;{\lie a})$ and $V(m_1,m_2;{\lie a})$ generate a subgroup of
finite index in $\SL(n,{\cal R})$.  
\end{proposition}
\begin{proof} Let $U$ be the sugroup generated by 
$H(m_1,m_2;{\lie a})$ and $V(m_1,m_2;{\lie a})$. Using 
the Steinberg relations (\ref{stein}) we conclude
that there is a non-zero two-sided ideal ${\lie b}$ ($={\lie a}^2$) in ${\cal
  R}$ such that every $E_{ij}(b)$ ($b\in {\lie b}$, $i\ne j$) is in $U$. The
main result of \cite{Vas} implies that $U$ has finite index in  
$\SL(n,{\cal R})$.
\end{proof}

\subsection{Constructing elements of ${\cal G}_{G,\pi}^1(\BZ)$}

Let $F_{n}$ ($n\in\BN,\, n\ge 3$)
be the free group with basis $x_1,\ldots ,x_{n-1},\, y$, $G$ a finite group
with generators
$g_1,\ldots , g_{n-1}$  and let
$$
\pi: F_{n}\to G\qquad \pi(x_1)=g_1,\ldots ,\pi(x_{n-1})=g_{n-1},\ \pi(y)=1
$$
be the redundant presentation resulting from the data. Notice the assumption
$n\ge 3$, if $n=2$ then $G$ is cyclic. We enclude this case into the next
section. 

Given $i=2,\ldots ,\ell$ and an ideal ${\lie a}$ of the order 
${\cal R}_i^\op$,
let us define 
\begin{equation}\label{pde1} 
H_i({\lie a}):=H_i(h_i,(n-2)h_i;{\lie a})\le \GL(n,\BZ)
\times \prod_{i=2}^\ell \GL(n-1,M(h_i,{\cal R}_i^\op))
\end{equation}
to be 
$H(h_i,(n-2)h_i;{\lie a})\le \GL((n-1)h_i,{\cal R}_i^\op)=
\GL(n-1,M(h_i,{\cal R}_i^\op))$ 
but considered as a subgroup of
the big direct product on the right hand side of (\ref{pde1}) 
(with all components of the factors different from the $i$-th factor equal to
the identity matrix). In case $i=1$ we have ${\cal O}_1=\BZ$ 
and we define 
$H_1({\lie a})=H_1(1,n-1;{\lie a})$ 
with $n-1$ replaced by $n$. Let us also
define 
\begin{equation}\label{pde2} 
V_i({\lie a}):=V_i(h_i,(n-2)h_i;{\lie a})\le \GL(n,\BZ)
\times \prod_{i=2}^\ell \GL(n-1,M(h_i,{\cal R}_i^\op))
\end{equation}
in a similar way. We have 
\begin{lemma}\label{eli1}
Let $i$ be one of the $1,\ldots ,\ell$ and $g\in H_i({\cal R}_i^\op)$. 
There is an
$e \in\BN$ such that $g^e$ is in $\sigma_{G,\pi}(\hat\Gamma(G,\pi))$
where $\sigma_{G,\pi}$ is as in (\ref{rep12}). 
\end{lemma}
\begin{proof} We prove the statement for $i\ge 2$, for $i=1$ the proof is only
different in notation.

Let $i\ge 2$ be fixed and consider $g\in H_i({\cal R}_i^\op)$. 
There are matrices $A_1,\ldots,A_{n-2}\in M(h_i,{\cal R}_i^\op)$ such that
the $i$-th component of $g_i$ of $g$ is equal to
$$g_i=\left(\begin{array}{cccccc}
           I_{h_i} &  0 & 0 & \ldots & 0 & 0\\
           0 & I_{h_i} & 0 & \ldots & 0 & 0\\
           \vdots &  \vdots &  \vdots &  \vdots &  \vdots &  \vdots  \\
           0 & 0 & 0 & \ldots & I_{h_i} & 0\\
           A_1 & A_2 & A_3 & \ldots & A_{n-2} &  I_{h_i}
\end{array} \right),
$$
all other components being equal to the identity matrix in the corresponding
component group.
By Corollary \ref{wichcor1}  
there are $B_1,\ldots,B_{n-1}\in \BQ[G]$ such that
$$\partial_1^{P_iB_1}+\ldots +\partial_{n-1}^{P_iB_{n-1}}({\rm Id}_{i,k})=A_k$$
for $k=1,\ldots, n-2$. Choose $e \in\BN$ to have sufficiently many divisors so
that $eP_iB_1,\ldots ,eP_iB_{n-2}$ are in $\BZ[G]$. Put 
$$\gamma:=\eta_{eP_iB_1,1}\circ \ldots \circ\eta_{eP_iB_{n-1},n-1} $$ 
with the $\eta$ defined in (\ref{eta}).
Increase now (if
necessary) $e$ to ensure $\gamma\in \hat\Gamma(G,\pi)$ (see (\ref{hut})). 
We have $\sigma_{G,\pi}(\gamma)=g^e$ and the lemma is proved.
\end{proof}
\begin{lemma}\label{eli2}
Let $i$ be one of the $1,\ldots ,\ell$ and $g\in V_i({\cal R}_i^\op)$. 
There is an
$e \in\BN$ such that $g^e$ is in $\sigma_{G,\pi}(\hat\Gamma(G,\pi))$.  
\end{lemma}
\begin{proof}
We prove the statement for $i\ge 2$, again for $i=1$ the proof is only
different in notation.

Let $i\ge 2$ be fixed and consider $g\in V_i({\cal R}_i^\op)$.
There are matrices $A_1,\ldots,A_{n-2}\in M(h_i,{\cal R}_i^\op)$ such that
the $i$-th component of $g_i$ of $g$ is equal to
$$g_i=\left(\begin{array}{cccccc}
           I_{h_i} &  0 & 0 & \ldots & 0 & A_1\\
           0 & I_{h_i} & 0 & \ldots & 0 & A_2\\
           \vdots &  \vdots &  \vdots &  \vdots &  \vdots &  \vdots  \\
           0 & 0 & 0 & \ldots & I_{h_i} & A_{n-2}\\
           0 & 0 & 0 & \ldots &  0 &  I_{h_i}
\end{array} \right),
$$
The matrices $A_1,\ldots,A_{n-2}$ can be considered as elements of 
$M(h_i,{\cal R}_i)$, after all $M(h_i,{\cal R}_i^\op)$ and 
$M(h_i,{\cal R}_i)$ are the same sets. Using the convention of Section
\ref{repsfini} we consider $u:=(A_1,\ldots,A_{n-2})$ as an element of 
$M(h_i,D_i)^{n-2}\subseteq \BQ\otimes_{\bb Z}\bar R_{\rm old}$ 
(see (\ref{gaschi})).   
Then there is  $e \in\BN$ such that $eu\in \bar R_{\rm old}$. 
We choose $U\in R_{\rm old} $ with $\bar U=eu$. We 
take $e \in\BN$ to have sufficiently many divisors so
that $\bar\gamma_U$ (see (\ref{ge})) is in $\hat\Gamma(G,\pi)$. We have 
$\bar\gamma_U=g^e$ and the lemma is proved.
\end{proof}

Since every $H_i({\cal R}_i^\op)$ or $V_i({\cal R}_i^\op)$ ($i=1,\ldots,\ell$)
is a finitely generated abelian group the two previous lemmas immediately 
imply
\begin{proposition}\label{prfi}
There is $e\in \BN$ such that
\begin{equation}\label{pde3}
\prod_{i=1}^\ell H_i(e{\cal R}_i)\le \sigma_{G,\pi}(\hat\Gamma(G,\pi))
\quad {and}\quad \prod_{i=1}^\ell V_i(e{\cal R}_i)\le 
\sigma_{G,\pi}(\hat\Gamma(G,\pi)).
\end{equation}
\end{proposition}

Finally we can state 
\begin{theorem}\label{prfit} Let $n\in \BN$ satisfy $n\ge 4$, then
$$\sigma_{G,\pi}(\hat\Gamma(G,\pi))\cap \SL(n,\BZ) 
\times \prod_{i=2}^\ell \SL((n-1)h_i,{\cal R}_i^\op))$$
is of finite index in $\SL(n,\BZ) 
\times \prod_{i=2}^\ell \SL((n-1)h_i,{\cal R}_i^\op))$.
\end{theorem}
\begin{proof} Just use Proposition \ref{prfi} together with Proposition 
\ref{wasser}.
\end{proof}

\section{Cyclic groups}\label{Choose}

In this section we present a detailed study of the case when 
$G$ is a cyclic group. We restrict ourselves to a particular 
type of (redundant) presentation, but following the discussion in Section
\ref{Concl3} this is in fact the general case. An important feature 
is that we are able to present a set of
generators of $\Gamma(G,\pi)$, in the case $n\ge 3$. For a general 
group $G$ we are not
able to do this. As a consequence we can determine the precise image of 
$\Gamma(G,\pi)$ (not only up to commensurability) 
acting on its invariant pieces of the relation module. 
Some of the results here will be needed in Section 
\ref{Proofs}, but we believe that this section is of independent interest.
In fact, our method can be observed here in a very simple, almost elementary
situation.   

Let us call the attention of the reader to the fact that our results
concerning cyclic groups, in fact, cover the general case of $G$ being
abelian. This follows since every irreducible representation of an abelian
group factors through a cyclic quotient, see the discussion in Section
\ref{SL3}.

Let the free group $F_n$ ($n\ge 2$) be generated 
by $x,\, y_1,\ldots, y_{n-1}$. Let us introduce the following elements of
$\A(F_n)$. Our convention here is again that values 
not given are identical to the argument.   
\begin{itemize}
\item $\delta_i(x):=y_ix$  and $\epsilon_i(x):=xy_i$ for $i=1,\ldots ,n-1$,
\item $\varphi_i(y_i):=xy_i$ and 
$\psi_i(y_i):=y_ix$ for $i=1,\ldots ,n-1$,
\item $\lambda_{ij}(y_i):=y_jy_i$ and $\nu_{ij}(y_i):=y_iy_j$ 
for $i,\, j=1,\ldots ,n-1$ with $i\ne j$.
\end{itemize}
A theorem of Nielsen (see \cite{MKS}) asserts that these elements generate 
$\A^+(F_n)$ where $\A^+(F_n)$ is the kernel of the homomorphism 
$\rho_1: \A(F_n) \to \GL(n,\BZ)$ followed by the determinant. If 
$\Gamma \le \A(F_n)$ is a subgroup, we define $\Gamma^+:=\Gamma\cap
\A^+(F_n)$. Let us further introduce 
\begin{itemize}
\item $\kappa_{jk}(x):= x [y_j,y_k]$, $\kappa_{ijk}(y_i):= y_i [y_j,y_k]$
for $1\le i,\, j,\, k\le n-1$ with $i\ne j$, $i\ne k$,
\item $\tau_{ij}(y_i):= y_i [x,y_j]$ for $1\le i,\, j\le n-1$ 
with $i\ne j$.
\end{itemize}
The set $T_n$ consisting of the $\kappa_{jk}$, $\kappa_{ijk}$, 
$\tau_{ij}$ (with indices as above), the 
$\delta_i^{-1}\circ \epsilon_i$, $\varphi_i^{-1}\circ \psi_i$ for
$i=1,\ldots ,n-1$ and the $\lambda_{ij}\circ\nu_{ij}^{-1}$ for
$1\le i,\, j\le n-1$ with $i\ne j$
generates the group ${\rm IA}(F_n)$ by
another theorem of Nielsen (see \cite{MKS}). Notice that $T_n$ is contained in
$\Gamma(G,\pi)$ for every presentation $\pi$ of an abelian group $G$.

Let us also introduce the following notation concerning matrix groups. We
write 
$$\Gamma^1(n,m)\le \SL(n,\BZ)$$
for the subgroup of $\SL(n,\BZ)$ consisting of those elements having a first
row which is congruent to $(1,0,\ldots,0)$ modulo $m\in\BN$. We need
\begin{lemma}\label{eli11} Let $n$ and $m$ be natural numbers with $n\ge
  3$. The group $\Gamma^1(n,m)$ is generated by the elementary matrices 
$$E_{ij}(1)\quad (i,\, j=1,\ldots,n, i\ne j, i\ne 1),\qquad E_{1j}(m) \quad
 (j=2,\ldots,n).$$   
\end{lemma}
This lemma is an elementary exercise using the euclidean algorithm in $\BZ$.

We first treat the cyclic group of order $2$ which we call $C_2$. We have
singled out this case since it is completely elementary and also plays a
special role in the proof of Corollary \ref{kaz}. 
Let $g$ be a generator of $C_2$. We choose
the presentation 
$$
\pi:F_n\to C_2,\quad \pi(x)=g,\ \pi(y_1)=\ldots=\pi(y_{n-1})=1.
$$
The corresponding relation module $\bar R$ has the following elements as a
$\BZ$-basis. 
$$\overline{x^2},\ \bar y_1,\ldots,\bar y_{n-1},\  
\overline{xy_1x^{-1}}=g\bar y_1,\ldots,
\overline{xy_{n-1}x^{-1}}=g\bar y_{n-1}.$$ 
The $\BQ$-vector space $\BQ\otimes_{\bb Z} \bar R$ decomposes as  
$\BQ\otimes_{\bb Z} \bar R=V_1\oplus V_{-1}$ where 
$V_1$, $V_{-1}$ are the $\pm 1$ 
eigenspaces of $g$ respectively.
Set $\bar R_1:=\bar R\cap V_1$ and $\bar R_{-1}:=\bar R\cap V_{-1}$. 
Notice that $\bar R_1+\bar R_{-1}$ is of finite index in $\bar R$. Introduce
\begin{equation}\label{basi}
v_i:=\bar y_i+\overline{xy_{i}x^{-1}},\quad 
w_i:=\bar y_i-\overline{xy_{i}x^{-1}}\qquad (i=1,\ldots ,n-1).
\end{equation}
Then $\overline{x^2},v_1,\ldots,v_{n-1}$ is a $\BZ$-basis of $\bar R_1$ and
$w_1,\ldots,w_{n-1}$ is a $\BZ$-basis of $\bar R_{-1}$. 
Since $\Gamma(C_2,\pi)$
leaves $\bar R_1$ and $\bar R_{-1}$ invariant we obtain, 
with the above $\BZ$-bases being chosen, representations
\begin{equation}\label{rep2}
\sigma_1:\Gamma^+(C_2,\pi)\to \GL(n,\BZ),\qquad 
\sigma_{-1}:\Gamma^+(C_2,\pi)\to \GL(n-1,\BZ).
\end{equation}
With all these data given we have: 
\begin{proposition}\label{cyc2} Let $n\ge 2$ be a natural number. The
  following hold.
 
{\rm (i)} The group $\Gamma^+(C_2,\pi)$ is generated by the 
automorphisms $\delta_i$, $\epsilon_i$
($i=1,\ldots,n-1$), the $\lambda_{ij}$ and $\nu_{ij}$ 
($i,\, j=1,\ldots ,n-1$, $i\ne j$), the $\varphi_i^2$ ($i=1,\ldots,n-1$) and
the elements of the set $T_n$ introduced above.

{\rm (ii)} The image of $\sigma_1$ is equal to $\Gamma^1(n,2)$.

{\rm (iii)} The index of $\Gamma^+(C_2,\pi)$ in $\A^+(F_n)$ is $2^{n}-1$.

{\rm (iv)} The representation $\sigma_{-1}$ is surjective onto $\GL(n-1,\BZ)$.

{\rm (v)} The group $\sigma_{-1}({\rm IA}(F_3))\le \GL(2,\BZ)$ is generated by
the matrices
$$\left(\begin{array}{cc}
           -1 &  0\\
            0  &  1
\end{array} \right),\qquad 
\left(\begin{array}{cc}
           1 &  0\\
           0 &  -1
\end{array} \right),\qquad
\left(\begin{array}{cc}
           1 &  2\\
           0 &  1
\end{array} \right),\qquad
\left(\begin{array}{cc}
           1 &  0\\
           2 &  1
\end{array} \right).
$$
\end{proposition}

The proof of this proposition is elementary but repeated below 
in the more general case, so we skip it here.
Items (ii), (iv) and (v) follow
 from formulas describing the action of the given
generators of $\Gamma^+(C_2,\pi)$ on the bases in $\bar R_1$ and $\bar R_{-1}$.
An important ingredient is Lemma \ref{eli11}. Note that $\Gamma^1(2,2)$ is
generated by the elementary matrices  
$$ E_{21}(1)=\left(\begin{array}{cc}
           1 &  0\\
           1 &  1
\end{array} \right),\qquad 
E_{12}(2)=\left(\begin{array}{cc}
           1 &  2\\
           0 &  1
\end{array} \right).$$
For $\Gamma^1(2,m)$ ($m\ge 3$) the analogous statement is not true.
Proposition \ref{cyc2} already
implies Corollary \ref{kaz}.  

We turn now to the case of a general cyclic group $C_m=\langle g\rangle$ 
($m\in \BN$). 
We choose the presentation 
$$
\pi:F_n\to C_m,\quad \pi(x)=g,\ \pi(y_1)=\ldots=\pi(y_{n-1})=1.
$$
The corresponding relation module $\bar R$ has the following 
$1+m(n-1)$ elements as a
$\BZ$-basis 
$$\overline{x^m},\   
\overline{x^{k}y_1x^{-k}},
\ldots, \overline{x^{k}y_{n-1}x^{-k}}\qquad (k=0,\ldots,m-1).$$
The rational group ring of $C_m$ decomposes as
\begin{equation}\label{grouri}
\BQ[C_m]=\prod_{d|m}\BQ(\zeta_d).
\end{equation}
Here $\zeta_d\in\BC$ is a primitive $d$-th root of unity. We consider
each of the $\BQ(\zeta_d)$ as a $C_m$-module letting $g$ act by multiplication
with $\zeta_d$. The ring of integers $\BZ(\zeta_d)$ is then a 
$\BZ[C_m]$-submodule of $\BQ(\zeta_d)$.
We mimic the decomposition (\ref{grouri}) inside the relation module 
by introducing
\begin{equation}
v_i(d,m):=\sum_{k=0}^{m-1}{\rm Tr}(\zeta_d^k)\, \overline{x^ky_ix^{-k}}
\end{equation}
for $i=1,\ldots,n-1$.
Here ${\rm Tr}(\zeta_d^k)$ is the trace of $\zeta_d^k$ taken from 
$\BQ(\zeta_d)$ to $\BQ$. Notice that the $v_i(d,m)$ are generalisations of the
generators in (\ref{basi}).
Let $\bar R_d$ be the $C_m$-submodule of $\bar R$ generated by the
$v_i(d,m)$ for $i=1,\ldots,n-1$ if $d>1$ and define 
$\bar R_1$ to be the $C_m$-submodule of $\bar R$ generated by the
$v_i(1,m)$ for $i=1,\ldots,n-1$ together with $\overline{x^m}$. By setting 
$$\Lambda_{d}(v_i(d,m)):=\left(0,\ldots,0,
\sum_{k=0}^{m-1}{\rm Tr}(\zeta_d^k)\zeta_d^k,0,\ldots,0\right)
\qquad (i=1,\ldots,n-1)$$
(with the non-zero entry in the $i$-th component) for $d>1$ and extending 
$C_m$-linearly we obtain well defined additive homomorphisms
$$\Lambda_d:\bar R_d\to \BZ(\zeta_d)^{n-1}\qquad (d|m,\, d>1). $$
The following are true 
\begin{itemize}
\item the homomorphism $\Lambda_d$ is $C_m$-equivariant and injective (onto an
  ideal in $\BZ(\zeta_d)$),
\item the $\BZ$-rank of $\bar R_d$ is $(n-1)\varphi(d)$ 
for $d>1$ and $n$ for $d=1$,
\item $\bar R_d\cap \bar R_{d'}=\langle 0\rangle$ for distinct divisors $d,\,
  d'$ of $m$,
\item the (direct) sum of all $\bar R_d$ ($d|m$) has finite index in $\bar R$,
\item  each $\bar R_d$ ($d|m$) is left invariant by $\Gamma(C_m,\pi)$.  
\end{itemize}
As above we transport the action of $\Gamma(C_m,\pi)$ on $\bar R_d$ 
to the image of $\lambda_d$ and obtain our representations
$$\sigma_{d}: \Gamma(G,\pi)\to \GL(n-1,\BQ(\zeta_d))\qquad (d|m,\, d>1).$$  
To control the image of $\sigma_{d}$ we first show:
\begin{proposition}\label{genem}
Let $n$ and $m$ be natural numbers with $n\ge 3$.  
The group $\Gamma^+(C_m,\pi)$ is generated by the 
automorphisms $\delta_i$, $\epsilon_i$
($i=1,\ldots,n-1$), the $\lambda_{ij}$ and $\nu_{ij}$ 
($i,\, j=1,\ldots ,n-1$, $i\ne j$) and the $\varphi_i^m$, 
($i=1,\ldots,n-1$) and
the elements of the set $T_n$ introduced above.
\end{proposition}
\begin{proof}
Since every element of $\Gamma(C_m,\pi)$ has 
to leave the normal closure of
$x^m$ in $F_n$ invariant 
and has to fix $x$ modulo $R$
we see that $\rho_1(\Gamma^+(C_m,\pi))$ has to be
contained in $\Gamma^1(n,m)$ (for $\rho_1$ see (\ref{furep})). 
Examining the action of elements in the
statement of the proposition on $F_n/F_n'$ and using Lemma \ref{eli11} we find
that $\rho_1(\Gamma^+(C_m,\pi))=\Gamma^1(n,m)$.

Now our set of generators contains enough elements to generate ${\rm IA}(F_n)$
by the theorem of Nielsen and also enough elements to generate the image of 
$\Gamma^+(G,\pi)$ in $\SL(n,\BZ)$ by Lemma \ref{eli11}. 
The proposition follows.
\end{proof}
Finally we can state:
\begin{proposition}\label{cycfin}
Let $n$ and $m$ be natural numbers with $n\ge 2$.
The image of the representation $\sigma_{d}$ ($d|m,\, d>1$) 
is $\GL^+(n-1,\BZ(\zeta_d))$
where $\GL^+(n-1,\BZ(\zeta_d))$ is the subgroup consisting of those elements in
$\GL(n-1,\BZ(\zeta_d))$ which have a power of $\zeta_d$
as determinant. 
\end{proposition}
\begin{proof}
Assume first that $n\ge 3$.
Evaluating $\sigma_{d}$ on the generators from Proposition
\ref{genem} we see that $\sigma_{d}(\Gamma^+(C_m,\pi))$ is contained in
$\GL^+(n-1,\BZ(\zeta_d))$.
We now show that all elementary matrices $E_{ij}(\zeta_d^k)$
($i,\, j=1,\ldots n-1,\, k\in \BZ$) are in the image of $\sigma_{d}$. This is
done by evaluating $\sigma_{d}$ on some of the generators from Proposition
\ref{genem}. We report only a special case: Take $\tau_{ij}$ for a pair
$(i,j)$ with $i\ne j$. We have
$\tau_{ij} (y_i)=y_i[x,y_j]$. This leads to
$$\bar\tau_{ij}(v_i(d,m))=v_i(d,m)+gv_j(d,m)-v_j(d,m)$$
which in turn shows that $E_{ij}(\zeta_d-1)$ is in the image of 
$\sigma_{d}$. Since $E_{ij}(1)$ is also present (use $\lambda_{ij}$) we have
shown that $E_{ij}(\zeta_d)$ is contained in the image of 
$\sigma_{d}$. 

After more elementary considerations like that we can show that all 
the elementary matrices with entries in $\BZ(\zeta_d)$ are in the image of 
$\sigma_{d}$. These matrices generate $\SL(n-1,\BZ(\zeta_d))$: For $n\ge 4$
this follows from results in \cite{Vas}, for $n=3$ and $d\ge 5$ we apply the
main result of \cite{Vas1}, in case $n=3$ and $d\le 4$ the ring $\BZ(\zeta_d)$
is euclidean, in case $n=2$ nothing has to be proved.

As a last step we use the $\varphi_i^{-1}\circ \psi_i$ ($i=1,\ldots, n-1$) to
show that the image of $\sigma_d$ contains a matrix with determinant equal to
$\zeta_d$. This finishes the proof for $n\ge 3$.

Let us finally discuss the case $n=2$ which seems trivial, but in fact is
not. Note that $\GL(1,\BZ(\zeta_d))$ contains elements of infinite order
whenever $d\ge 5$, $d\ne 6$ holds. Since we do not have generators for 
$\Gamma^+(2,m)$ ($m\ge 3$) at hand we cannot directly rule out the 
possibility that such an element lies in the image of $\sigma_d$. We argue as
follows. Let $\varphi\in \Gamma^+(2,m)$ be such that for some divisor $d$ of
$m$ the image $\sigma_d(\varphi)$ has a determinant of infinite order. We
extend $\varphi$ to an element of $\psi$ of $\Gamma^+(3,m)$ by setting
$\psi(x):=\varphi(x)$, $\psi(y_1):=\varphi(y_1)$, $\psi(y_2):=y_2$. Since the
determinants of $\sigma_d(\varphi)$ and $\sigma_d(\psi)$ coincide we have
finished also the case $n=2$.  
\end{proof}
Proposition \ref{cycfin} allows us to conclude the following result concerning
the abelianisations of the images of $\Gamma(C_m,\pi)$ under the
representations $\sigma_d$ ($d|m,\, d>1$).  
\begin{corollary}\label{cycfincor}
Let $n$ and $m$ be natural numbers with $n\ge 2$.
Then $\sigma_{d}(\Gamma(C_m,\pi))$ ($d|m,\, d>1$) has finite abelianisation.
\end{corollary}
\begin{proof}
It is well known that $\GL^+(n-1,\BZ(\zeta_d))$ has finite abelianisation for
all $d\in\BN$: For $n=2$ this fact is obvious,
for $n=3$ and $d=2,\, 3,\, 4,\, 6$ a presentation of 
$\SL(2,\BZ(\zeta_d))$ (see \cite{EGM}) can be used, for $d=5$ or $d\ge 7$ this
fact is contained in \cite{Serresl}. 
For other values of $n,\, d$ the group $\GL^+(n-1,\BZ(\zeta_d))$ has Kazdhan's
property (T).  
Since $\sigma_{d}(\Gamma(C_m,\pi))$ contains 
$\sigma_{d}(\Gamma^+(C_m,\pi))$ as a subgroup of finte index the result
follows.   \end{proof}

We finally remark that
$$\sigma_{d} ({\rm IA}(F_n))\le \GL^+(n-1,\BZ(\zeta_d))\qquad (d|m,\, d>1) 
$$
is equal to the full congruence subgroup of $\GL^+(n-1,\BZ(\zeta_d))$ with
respect to the ideal of $\BZ(\zeta_d)$ generated by $\zeta_d-1$. This can be
shown by a more detailed analysis of the arguments involved in the proof of
Proposition \ref{cycfin}.

\section{Completion of proofs}\label{Proofs}

In this section we complete the proofs of all the results promised in the
introduction.

Theorem \ref{prfit} proves Theorem \ref{teo}. Both theorems were formulated
only for $n\ge 4$ but our method applies also to the
case $n=3$. The only obstacle for $n=3$ is that congruence elementary matrices
in $\SL(2)=\SL(n-1)$ do not always generate a finite index subgroup, for
example the matrices 
$$ \left(\begin{array}{cc}
           1 &  3\\
           0 &  1
\end{array} \right),\qquad 
\left(\begin{array}{cc}
           1 &  0\\
           3 &  1
\end{array} \right)$$
do not generate a subgroup of finite index in $\SL(2,\BZ)$. By the result of
Vaserstein \cite{Vas1} this is a real obstacle (in the case of commutative
rings of endomorphisms) only for $\SL(2,\BZ)$ and $\SL(2)$ of a ring of
integers in an imaginary quadratic number field. Thus our results still work
for $n=3$ in many cases. For example
if the 
group algebra $\BQ[G]$ of the finite group $G$ has 
(except for the trivial module) no irreducible module $N$ with
${\rm dim}_D(N)=1$ ($D={\rm End}_G(N)$). If there are such  modules $N$ but
if in each case $D$ is a number field which is not $\BQ$ 
or imaginary quadratic we are also fine. The other cases remain open.

\smallskip
\begin{prof} {\it of Theorem \ref{teo4}:} $n=2$: In this case we 
just consider the (surjective) representation
$$\rho_1 : \A(F_2)\to \A(F_2/F'_2)\cong \GL(2,\BZ).$$
It is well known that $\GL(2,\BZ)$ contains a free subgroup 
of any rank ($\ge 2$) of finite index. Hence any finitely generated group
is the image of a subgroup of finite index in $\A(F_2)$. Note that the
groups appearing in Theorem \ref{teo4} as image groups are finitely generated
since they are arithmetic groups.   

$n=3$: Here we use the representation 
$\sigma_{-1}:\Gamma^+(C_2,\pi)\to \GL(2,\BZ)$ 
from Section \ref{Choose} to find 
for every $r\in \BN$ with $r\ge 2$ a subgroup of finite index in $\A(F_3)$
which can be mapped onto the free group of rank $r$. We then argue as in the
case $n=2$.

$n\ge 4$: We first consider the case $k=1$. For $h,\, m\in\BN$ we have to find
a subgroup $\Gamma\le \A(F_n)$ of finite index and a representation 
$$
\rho : \Gamma \to \SL((n-1)h,\BQ(\zeta_{m}))^{m} 
$$
such that $\rho(\Gamma)$ is commensurable with 
$\SL((n-1)h,\BZ(\zeta_{m}))^{m}$. To do this we take the group
$$G:=S_{h+1}\times C_m\times C_m.$$ 
As an input for Theorem \ref{prfit} we need  suitable $\BQ[G]$-modules. We let
$PM_{h+1}$ be the standard $h+1$ dimensional $\BQ[S_{h+1}]$-permutation
module and let $M_h$ be the kernel of the augmentation map from 
$PM_{h+1}$ to $\BQ$. The $\BQ[S_{h+1}]$-module $M_h$ is 
irreducible and has $\BQ$-dimension
$h$. We put $N:=\BQ(\zeta_m)\otimes_\BQ M_h$ and consider $N$ as a 
$\BQ$-vector space. Here, as before, $\zeta_m\in\BC$ is a primitive $m$-th
root of unity. Let $\epsilon : C_m\times C_m\to \langle \zeta_m\rangle$ be a
surjective homomorphism. We turn the $\BQ$-vector space 
$N=\BQ(\zeta_m)\otimes_\BQ M_h$ into a $\BQ[G]$-module $N_\epsilon$ by
letting $S_{h+1}$ act on $M_h$ and letting $g\in C_m\times C_m$ act by
multiplication by $\epsilon(g)$ on $\BQ(\zeta_m)$. The following are clear
for every surjective homomorphism 
$\epsilon : C_m\times C_m\to \langle \zeta_m\rangle$.
\begin{itemize}     
\item ${\rm dim}_\BQ(N_\epsilon)=\varphi(m)h$,
\item ${\rm End}_G(N_\epsilon)=\BQ(\zeta_m)$,
\item ${\rm dim}_{\BQ(\zeta_m)}(N_\epsilon)=h$,
\item if $\epsilon_1,\,\epsilon_2:  C_m\times C_m\to \langle \zeta_m\rangle$ 
are two surjective homomorphisms with distinct kernels then $N_{\epsilon_1}$
and $N_{\epsilon_2}$ are not isomorphic as $\BQ[G]$-modules. 
\end{itemize}     
Our group $G$ can be generated by $3$ elements: take $g_1,\, g_2,\, g_3\in 
S_{h+1}$ such that $g_3$ and the commutator of $g_1,\, g_2$ generate 
$S_{h+1}$ and consider the elements 
$(g_1,1,0)$, $(g_2,0,1)$, $(g_3,0,0)$ in $G$. Let $\pi :F_n\to G$ be a
redundant presentation of $G$ (remember the assumption is $n\ge 4$). 
Since the number of surjective homomorphisms 
$\epsilon : C_m\times C_m\to \langle \zeta_m\rangle$ with pairwise distinct
kernels is bigger or equal to $m$ we may apply Theorem \ref{prfit} to obtain
a subgroup $\Gamma \le \A(F_n)$ and a representation 
$$\sigma :\Gamma\to   \GL((n-1)h,\BQ(\zeta_m))^m$$
such that $\sigma(\Gamma)\cap \SL((n-1)h,\BZ(\zeta_m))^m$ is of finite index 
in $\SL((n-1)h,\BZ(\zeta_m))^m$. Here we do not know whether 
$\sigma(\Gamma)\cap \SL((n-1)h,\BZ(\zeta_m))^m$ is of finite index in   
$\sigma(\Gamma)$. However,
there is a subgroup $\Delta\le \GL((n-1)h,\BZ(\zeta_m))^m$ of finite index
such that $ \SL((n-1)h,\BZ(\zeta_m))^m\le \Delta$ and such that
$\Delta\le Z\cdot \SL((n-1)h,\BZ(\zeta_m))^m$ with a subgroup 
$Z\le \GL((n-1)h,\BZ(\zeta_m))^m$ which is central and satisfies 
$Z\cap \SL((n-1)h,\BZ(\zeta_m))^m=\langle 1\rangle$. Let $\Gamma_0\le \Gamma$
be the inverse image of $\Delta$ under $\sigma$. The representation 
$$\sigma_0:  \Gamma_0\to   \Delta/Z=\SL((n-1)h,\BZ(\zeta_m))^m$$
satisfies the requirements for $k=1$. Another way of looking at the last
construction is to project a finite index subgroup of the image of 
$\sigma(\Gamma)$ into ${\rm PGL}(n-1,\BZ(\zeta_m))$ which contains 
${\rm PSL}(n-1,\BZ(\zeta_m))$ as a finite index subgroup. We then intersect
this subgroup in ${\rm PGL}(n-1,\BZ(\zeta_m))$ with 
${\rm PSL}(n-1,\BZ(\zeta_m))$ and pull the resulting subgroup of 
${\rm PSL}(n-1,\BZ(\zeta_m))$ back to obtain a finite index subgroup of
$\Gamma$ which is mapped into ${\rm SL}(n-1,\BZ(\zeta_m))$.

We turn now to the general case. Let $\Gamma_1,\ldots,\Gamma_k\le \A(F_n)$
be the subgroups of finite index and 
$$\sigma_i:  \Gamma_i\to   \SL((n-1)h_i,\BZ(\zeta_{m_i}))^{m_i}\qquad
(i=1,\ldots,k)$$
be the representations constructed above. Set 
$\Gamma:=\Gamma_1\cap\ldots\cap\Gamma_k$ and consider the representation
$$\sigma:  \Gamma\to   \prod_{i=1}^k\SL((n-1)h_i,\BZ(\zeta_{m_i}))^{m_i}$$
which maps an element $\gamma\in\Gamma$ to the tuple
$(\sigma_1(\gamma),\ldots,\sigma_k(\gamma))$. The image of $\Gamma$ projects
onto an arithmetic subgroup in each of the factors of the direct product.
We finish the proof by applying the following lemma.
\begin{lemma}\label{Margu}
Let ${\cal H}_1,\ldots ,{\cal H}_k$ be pairwise non-isogeneous
$\BQ$-defined simple linear algebraic groups which all have real rank greater
or equal to $2$. Let
$$\Gamma\le \prod_{i=1}^k{\cal H}_i^{m_i}(\BQ)$$ 
be a subgroup such that all projections into the factors ${\cal
  H}_i^{m_i}(\BQ)$ are arithmetic groups. Then $\Gamma$ is an arithmetic
subgroup of $\prod_{i=1}^k{\cal H}_i^{m_i}$, that is $\Gamma$ is commensurable
with $\prod_{i=1}^k{\cal H}_i^{m_i}(\BZ)$.
\end{lemma}
This lemma is proved using Margulis-superrigidity. See \cite{L} item (iv) on
pages 328-329 or
\cite{GP} Section 2. 
\end{prof}

\section{(SL) or not (SL)?}\label{SL}

In this section we show some results concerning
the question whether $\rho_{G,\pi}(\Gamma(G,\pi))\cap{\cal G}^1_{G,\pi}$ is of
finite index in $\rho_{G,\pi}(\Gamma(G,\pi))$ or not (for the notation see the
introduction). The relevance of this 
question and the notation are explained in the introduction. We change 
here from
the consideration of the rational relation module to its complex version. 
This brings certain technical advantages and still 
reflects the problem.  

Let us set up some notation for this section. Let $G$ be a finite group and
$\pi:F_n\to G$ an epimorphism with kernel $R$. We define $\Gamma(G,\pi)$, 
${\cal G}_{G,\pi}$ and ${\cal G}^1_{G,\pi}$ as in the introduction and
consider the action of $\Gamma(G,\pi)$ now on the complex relation module
$\BC\otimes_\BZ\bar R$. The result of Gasch\"utz now says that
$$\BC\otimes_{\bb Z}\bar R=\BC\oplus \BC[G]^{n-1}$$ 
as $\BC[G]$-modules.

If $Q$ is an irreducible $\BC[G]$-module and 
$M$ is an arbitrary $\BC[G]$-module we keep the notation ${\bf I}_Q(M)$ 
for the $Q$-isotypic component inside $M$. We hope that this creates no
confusion with the ${\bf I}_N(M)$ of Section \ref{sura} where $N,\, M$ are 
$\BQ[G]$-modules.

Let $Q$ be an non-trivial irreducible $\BC[G]$-module of complex dimension 
$h_Q$. Let ${\cal G}_{G,\pi,Q}$ be the stabiliser of 
$I_Q(\BC\otimes_\BZ\bar R)$ inside ${\cal G}_{G,\pi}$ (see \ref{stabi}).
By our usual procedure we obtain a representation 
\begin{equation}
\rho_{G,\pi,Q}:\Gamma(G,\pi)\to {\cal G}_{G,\pi,Q}=\GL(n-1,M(h_Q,\BC))=
\GL((n-1)h_Q,\BC).
\end{equation}
We call the pair $(G,Q)$ of {\it type (SL)} if 
$\rho_{G,\pi,Q}(\Gamma(G,\pi))\cap  \SL((n-1)h_Q,\BC)$ is of finite index in 
$\rho_{G,\pi,Q}(\Gamma(G,\pi))$ for every $n\in \BN$ and every epimorphism
$\pi:F_n\to G$. We say the finite group $G$ is of {\it type (SL)}
if for all non-trivial irreducible $\BC[G]$-module $Q$ the pair 
$(G,Q)$ is of type (SL). 

Note that, since $\Gamma(G,\pi)$ is finitely generated, 
$\rho_{G,\pi,Q}(\Gamma(G,\pi))\cap  \SL((n-1)h_Q,\BC)$ is of finite index in 
$\rho_{G,\pi,Q}(\Gamma(G,\pi))$ if and only if $\rho_{G,\pi,Q}(\Gamma(G,\pi))$
contains no element with a determinant of infinite order.

Note that the above implies that once we find a finite group $G$
which is not of type (SL), we have found a subgroup $\Gamma$ of finite 
index in some $\A(F_n)$ which has an epimorphism
onto $\BZ$.

Let $G$ be a finite group $G$ and $N_2,\ldots,N_\ell$ its non-trivial
irreducible $\BQ[G]$-modules. Let $\pi:F_n\to G$ be a surjective homomorphism.
The following are equivalent: 
\begin{itemize}
\item For all $i=2,\ldots,\ell$ the intersection  
$$\sigma_{G,\pi,i}(\Gamma(G,\pi))\cap \SL((n-1)h_i,{\cal R}_i^{\rm op})$$
has finite index in  
$\sigma_{G,\pi,i}(\Gamma(G,\pi))\le \GL((n-1)h_i,{\cal R}_i^{\rm op})$.
\item For all non-trivial  
irreducible $\BC[G]$-modules $Q$ the intersection 
$$\rho_{G,\pi,Q}(\Gamma(G,\pi))\cap \SL((n-1)h_Q,\BC)$$
 is of finite index in $\rho_{G,\pi,Q}(\Gamma(G,\pi))$.
\end{itemize}
This is seen by introducing the decomposition of the
complexifications $\BC\otimes_\BQ N_i$ into irreducible $\BC[G]$-modules.  

There are some non-trivial irreducible $\BC[G]$-modules $Q$ which are
obviously of type (SL):
\begin{proposition}\label{aal}
Let $G$ be a finite group and $Q$ a non-trivial irreducible $\BC[G]$-module.
If there is an (irreducible) $\BQ[G]$-module $N$ with $Q=\BC\otimes_\BQ N$
then $(G,Q)$ is of type (SL).
\end{proposition}
\begin{proof}
We note that under the above assumptions the endomorphism algebra 
$D:={\rm End}_G(N)$ is a division algebra 
over the rational numbers with $\BQ$ as
  center. In this case $\SL(m,{\cal R})$ (i.e. the group of elements of
  reduced norm $1$ in $M(m,{\cal R})$) is of finite index in 
$\SL(m,{\cal R})$ for any order ${\cal R}\le D$ and any $m\in \BN$.
\end{proof}
In case of the symmetric groups $S_m$ ($m\in \BN$) it is known, that every
$\BC[S_m]$-module satisfies the assumptions of Proposition \ref{aal}, see
\cite{S}, Section 13.  This implies:
\begin{corollary}
All the symmetric groups $S_m$ are of type (SL).
\end{corollary}
In case of the alternating group $A_5$ Proposition \ref{aal} applies to
the irreducible $\BC[A_5]$-modules of dimensions $4$ and $6$ but not to the
two $3$-dimensional irreducible $\BC[A_5]$-modules. They fit together to form
a $6$-dimensional irreducible $\BQ[A_5]$-module $N$ with 
${\rm End}_G(N)=\BQ(\sqrt{5})$.

\subsection{Reduction Steps}\label{SL2}

After the definitions above we discuss some reduction steps which will allow 
us to recognise pairs $(G,Q)$ of type (SL).

\begin{lemma}\label{redu01}
Let $G$ be a finite group with normal subgroup $H_0$ and quotient $H=G/H_0$. 
Let $Q$ be a 
non-trivial irreducible $\BC[G]$-module on which $H_0$ acts trivially. If 
$(H,Q)$ is of type (SL) then $(G,Q)$ also has this property.
\end{lemma}
\begin{proof}  Let $\pi:F_n\to G$ ($n\ge 2$) be any epimorphism 
and let $R_G$ be the kernel of $\pi$. Let $\pi_0 :F_n\to G\to H$ be the
resulting presentation of $H$ and $R_H$ its kernel. Obviously $R_G$ is a
subgroup of finite index in $R_H$. It is also easily seen that $\Gamma(G,\pi)$
is a subgroup of finite index of $\Gamma(H,\pi_0)$. The inclusion of   
$R_G$ into $R_H$ gives rise to a surjective $\BC$-linear map
$$\Theta: \BC\otimes_{\bb Z} \bar R_G\to \BC\otimes_{\bb Z} \bar R_H.$$
This linear map is equivariant for the induced actions of $G$ and 
$\Gamma(G,\pi)$.
By restriction we obtain a map
\begin{equation}\label{fro1}
\Theta_Q: {\bf I}_Q(\BC\otimes_{\bb Z} R_G)\to{\bf I}_Q(\BC\otimes_{\bb Z} R_H)
\end{equation}
which is surjective and $H$-equivariant. 
Gasch\"utz' result implies that $\Theta_Q$ is in
fact an ($\Gamma(G,\pi)$-equivariant) isomorphism.  

We have to show that ${\rm det}(\rho_{G,\pi,Q}(\varphi))$ 
is of finite order
for any $\varphi\in\Gamma(G,\pi)$. We compute this determinant on the right
hand side of (\ref{fro1}) where our assumption implies this property.
\end{proof}

\begin{lemma}\label{redu2}
Let $G$ be a finite group with  normal subgroup $H$. Let $Q$ be a 
non-trivial irreducible $\BC[G]$-module which is induced from $H$, that is,
$Q$ is of the form $Q={\rm Ind}_H^G(Q_0)$ where $Q_0$ is an irreducible
$\BC[H]$-module. If 
$(H,Q_0)$ is of type (SL) then $(G,Q)$ also has this property.
\end{lemma}
\begin{proof} Let $\pi:F_n\to G$ ($n\ge 2$) be any epimorphism 
and let $R$ be the kernel of $\pi$. The subgroup 
$\pi^{-1}(H)\le F_n$ is a free group on $m=1+[G:H](n-1)$ generators, 
we denote it
by $F_m$. We then have a commutative diagram  
$$\xymatrix{
1 \ar[r] & R  \ar[r]  & F_n \ar[r]^{\pi} & G \ar[r] & 1\\                      
1 \ar[r] & R  \ar[u]^{=}\ar[r] & F_m \ar[u]\ar[r]^{\pi_0} & H\ar[u]\ar[r] & 1\\
}
$$
where $\pi_0$ is the restriction of $\pi$ to $F_m$. Note that 
$R={\rm ker}(\pi)={\rm ker}(\pi_0)$.  

An element $\varphi\in\Gamma(G,\pi)$ stabilises the subgroup $F_m\le F_n$ and
as an automorphism of $F_m$ it lies in $\Gamma(H,\pi_0)$. On the vector space 
$\BC\otimes_\BZ\bar R$ the linear maps induced by $\varphi$ as an element of 
$\Gamma(G,\pi)$ and by $\varphi$ as an element of $\Gamma(H,\pi_0)$ coincide.
We denote both by $\bar\varphi$. 

The complex relation module $\BC\times_\BZ\bar R$ now has a structure as
a $\BC[G]$-module and a structure of a $\BC[H]$-module. In fact the second is
just the restriction of the first. 

Let $Q=Q_0\oplus Q_1\oplus\ldots \oplus Q_r$ be the decomposition of $Q$ into
irreducible $\BC[H]$-modules 
We now prove that 
\begin{equation}\label{frob}
{\bf I}_Q(\BC\otimes_{\bb Z}\bar R)=
{\bf I}_{Q_0}(\BC\otimes_{\bb Z}\bar R)\oplus
{\bf I}_{Q_1}(\BC\otimes_{\bb Z}\bar R)\oplus\ldots \oplus 
{\bf I}_{Q_r}(\BC\otimes_{\bb Z}\bar R).
\end{equation}
On the left hand side the isotypic component of the $\BC[G]$-module $Q$ is
taken whereas on the right hand side we see the the isotypic component of the 
corresponding $\BC[H]$-modules. Clearly the left hand side of (\ref{frob})
is contained in the right hand side.
Note that, since $H$ is normal in $G$ the module $Q$ is also induced from any
of the $Q_0,\ldots, Q_r$. 
The Frobenius reciprocity law (\cite{S}, chapter 7) 
implies that none of the $Q_0,\ldots, Q_r$ can occur as an irreducible
constituent of the restriction of an irreducible $\BC[G]$-module different
from $Q$. 

We have to show that ${\rm det}(\rho_{G,\pi,Q}(\varphi))$ is of finite order
for any $\varphi\in\Gamma(G,\pi)$. We compute this determinant on the right
hand side of (\ref{frob}) where our assumption implies this property.
\end{proof}

\subsection{Abelian and metabelian groups}\label{SL3}

This section contains results which describe the image of $\rho_{G,\pi}$ for
abelian and then for metabelian groups $G$ up to commensurability. In Section  
\ref{Choose} we have proved a very precise result in this direction for cyclic
groups $G$. This allowed us to conclude
for $n\ge 2$ that the image of $\Gamma(G,\pi)$ 
(under certain linear representations) has finite abelianisation. 
We present here the following much stronger result for $n\ge 3$. This result
is due to A. Rapinchuk (it appears in a somewhat hidden form in \cite{RA},
page 150), we thank him for the permission to include it. 
\begin{proposition}\label{rapi}
Let $n\ge 3$ be a natural number, $G$ a finite abelian group $G$ and let
$\pi:F_n\to G$ be any surjective homomorphism. Then $\Gamma(G,\pi)$ has finite
abelianisation. 
\end{proposition}
\begin{proof}
Since $G$ is abelian we have $F_n'\le R$ where $R$ is, as always, the kernel
of $\pi$. Consequently we have ${\rm IA}(F_n)\le \Gamma(G,\pi)$. Writing
$\tilde\Gamma(G,\pi):=\rho_1(\Gamma(G,\pi))\le \GL(n,\BZ)$ for the image of 
$\rho_1(\Gamma(G,\pi))$ in $\GL(n,\BZ)$ we have the exact sequence of groups
$$1\to {\rm IA}(F_n)\to \Gamma(G,\pi)\to \tilde\Gamma(G,\pi)\to 1$$
Note that $\tilde\Gamma(G,\pi)$ is a subgroup of finite index in 
$\GL(n,\BZ)$. Since $\Gamma(G,\pi)$ is a finitely generated group we need only
show that 
$$H^1(\Gamma(G,\pi),\BQ)={\rm Hom}(\Gamma(G,\pi),\BQ)=0$$
We shall do that by applying the edge-term sequence of the Hochschild-Serre
spectral sequence for the above short exact sequence of groups. We shall use
that there is an exact sequence of $\BQ$-vector spaces
\begin{equation}\label{Hose}
H^1(\tilde\Gamma(G,\pi),\BQ)\to H^1(\Gamma(G,\pi),\BQ)\to 
H^1({\rm IA}(F_n),\BQ)^{\tilde\Gamma(G,\pi)}
\end{equation}
First of all the left hand side of (\ref{Hose}) is $0$ since $n\ge 3$ and
$\GL(n,\BZ)$ has Kazdhan's property (T).

Next we evaluate the right hand side of  (\ref{Hose}). By a theorem of
Bachmuth (\cite{Ba} and also \cite{F}) there is a $\A(F_n)$-equivariant
isomorphism 
$$H^1({\rm IA}(F_n),\BZ)\to {\rm Hom}_{\bb Z}(F_n/F_n',F_n'/\gamma_2(F_n))$$
where $\gamma_2(F_n):=[[F_n,F_n],F_n]$ is the third term of the lower central
series of $F_n$. Let $\GL^\pm (n,\BC)$ be the subgroup of   
$\GL (n,\BC)$ consisting of those elements having determinant $\pm 1$. Let
$\BC^n$ be the standard $\GL^\pm (n,\BC)$-module of dimension $n$ and
$\Lambda^2(\BC^n)$ be its outer square. Consider $\GL^\pm (n,\BZ)$ as a
subgroup of $\GL^\pm (n,\BC)$. By identifying $F_n'/\gamma_2(F_n))$ with
$\Lambda^2(F_n/F_n')=\Lambda^2(\BZ^n)$ we obtain an
$\GL^\pm (n,\BZ)$-equivariant
isomorphism  
$$\BC\otimes_{\bb Z} {\rm Hom}_{\bb Z}(F_n/F_n',F_n'/\gamma_2(F_n))\to
{\rm Hom}(\BC^n,\Lambda^2(\BC^n)).$$
It is well known (see \cite{F}, page 427)
that the right hand side splits as a $\GL^\pm (n,\BC)$-module
as 
$${\rm Hom}(\BC^n,\Lambda^2(\BC^n))=\BC^n\oplus V_n$$
where $V_n$ is a $\GL^\pm (n,\BC)$-module of dimension 
${\rm dim}_\BC(V_n)=n(n+1)(n-2)/2$ which is irreducible even as 
$\SL(n,\BC)$-module. Hence the space of $\SL(n,\BC)$-invariants in   
${\rm Hom}(\BC^n,\Lambda^2(\BC^n))$ is $0$. Since the Zariski closure of
$\tilde\Gamma(G,\pi)$ contains $\SL(n,\BC)$ the group 
${\rm Hom}_\BZ(F_n/F_n',F_n'/\gamma_2(F_n))^{\tilde\Gamma(G,\pi)}$ is $0$.
This shows that the right hand side of (\ref{Hose}) is $0$. Together with the
above the proposition is proved.
\end{proof}
Proposition \ref{rapi} is clearly not true for $n=2$, since the left-most term
in (\ref{Hose}) contributes to the $H^1(\Gamma (G,\pi),\BQ)$. See also the
tables in Section \ref{Comp}. 
As an immediate corollary from Proposition \ref{rapi} we get:
\begin{corollary}
All finite abelian groups are of type (SL).
\end{corollary}
This corollary can also be deduced from Corollary \ref{cycfincor} (in fact
this corollary is needed in the case $n=2$) and  Lemma \ref{redu01}.
But the statement
incorporated in Proposition \ref{rapi} (for $n\ge 3$) is in fact stronger.

A metabelian group is a group $G$ with abelian commutator subgroup $G'$.
We prove:
\begin{proposition}
All metabelian finite groups $G$ are of type (SL)
\end{proposition}
\begin{proof}
We use induction on the order of $G$ and the fact proved above that all finite
abelian groups have property (SL).

Let $Q$ be a non-trivial irreducible $\BC[G]$-module.
We may assume by induction and Lemma \ref{redu01} that $G$ acts faithfully on
$Q$. 

{Case 1:} the commutator subgroup is not contained in the center ${\rm Z}(G)$ 
of $G$. The restriction of the module $Q$ to $G'$ cannot be
isotypic. Otherwise the abelian subgroup $G'$ would act by scalar matrices on
$Q$ and $G'$ would have to be contained in the center of $G$. 
Now Proposition 24 of \cite{S} says that $Q$ is either isotypic or induced.
It follows that
$Q$ is induced from a proper subgroup $H$ of $G$ containing $G'$.

By our construction $H$ is normal in $G$. 
Hence this case is finished by induction and
application of Lemma \ref{redu2}.   

{Case 2:} the commutator subgroup is contained in the center ${\rm Z}(G)$ 
of $G$. In this case $G$ is nilpotent of class $2$. Unless $G$ is abelian (the
case treated above) we find $g\in G$ with $g\notin {\rm Z}(G)$. The subgroup
$A:=\langle G',g\rangle$ generated by $G'$ and $G$ is abelian, normal in $G$
and not contained in ${\rm Z}(G)$. We then proceed as in case 1. 
\end{proof}

\section{Concluding remarks and suggestions for further research}\label{Concl}

\subsection{Dependence on the presentation}\label{Concl3}

The results of the current paper can be considered as a first step toward a
systematic study of the relation module of a finite group as a
$\Gamma(G,\pi)$-module, i.e. an equivariant Gasch\"utz's theory. Here we gave
a quite satisfactory answer in the case of redundant presentations. It is not
clear to what extend this represents the general case. 

For a fixed finite group $G$ let us look at all the possible epimorphisms from
the free group $F_n$ to $G$ and ${\bf R}(n,G)$ the set of their kernels. The
automorphism group $\A(F_n)$ acts on this set. It is easy to see that the
equivariant Gasch\"utz theory that we are trying to develop here depends only
on the orbits of  $\A(F_n)$ on ${\bf R}(n,G)$ and not on the actual
presentation. Various authors studied the transitivity properties 
of the action of $\A(F_n)$ on ${\bf R}(n,G)$. It is not transitive in general 
(a first example was given by B. H. Neumann \cite{Neu})
but it is in some interesting cases. For example, it is transitive if $G$ is a
cyclic group and $n\ge 2$. Thus our results in Section \ref{Choose} give the
full picture for all relation modules of cyclic groups. It is also
transitive if $n>2{\rm log}_2(|G|)$ or if $G$ is solvable and $n>d(G)$, where
$d(G)$ denotes the minimal number of generators of $G$ (see \cite{DU}). An
old conjecture of Wiegold predicts transitivity for finite simple groups if
$n\ge 3$, in which case again $n>d(G)$ holds. For some partial results see
\cite{Gi}, \cite{Ev} and the references in \cite{Pak}. 
It was even suggested in \cite{Pak}
that for every finite group $G$ it is transitive if $n>d(G)$. 

In \cite{Ev} it is shown that if a finite group $G$ 
has spread 2 (i.e. for every
$g_1,\, g_2\in G$ with $g_1\ne 1 \ne g_2$
there is  $h\in G$ such that $\langle g_1,\, h\rangle=    
\langle g_2,\, h\rangle=G$) then $\A(F_n)$ acts transitively on the set of
kernels of all redundant presentations. Guralnick and Shalev  
\cite{GS} proved that almost all
finite simple groups have spread 2.

Theorem 1.4 is not true in general if the presentation is not redundant, at
least when $n=2$, see the $A_5$-example in Section \ref{Comp}.
The story for $n\ge 3$ may be different as $\A(F_2)$ behaves
in many ways differently from $\A(F_n)$ for $n\ge 3$.

\subsection{Representations of IA($F_n$)}\label{Concl1}

A well known theorem of Formanek and Procesi (see \cite{FP}) asserts 
that the automorphism group $\A(F_n)$ has, for $n\ge 3$, no faithful
linear representation over any field. Their proof suggests a stronger  

{\it Conjecture: For every linear representation $\rho$ of $\A(F_n)$ the image 
$\rho({\rm Inn}(F_n))$ is virtually solvable.}

As far as we know, in all previously constructed cases of linear
representations even 
$\rho({\rm IA}(F_n))$ had this property. 
It should be mentioned that in all the representations $\rho$ studied in this
paper, $\rho({\rm Inn}(F_n))$ is finite but $\rho({\rm IA}(F_n))$ is far from
being virtually solvable and, in fact, our results show that 
$\rho({\rm IA}(F_n))$ can be Zariski dense in very large semi-simple groups. 
Our reprentations seem to be the first known with this property. 
We have:   
\begin{proposition}\label{teo3}
Let $n\ge 3$, $k\ge 1$, $h_1<\ldots <h_k$, $m_1,\ldots ,m_k$ be 
natural numbers.  Let $\BQ(\zeta_{m_i})$ be the field of 
$m_i$-th roots of unity and
$\BZ(\zeta_{m_i})$ its ring of integers. There is a 
subgroup $\Gamma\le {\rm IA}(F_n)$ of finite index and a representation 
$$
\rho : \Gamma \to \prod_{i=1}^k\SL((n-1)h_i,\BQ(\zeta_{m_i}))^{m_i} 
$$
such that $\rho(\Gamma)$ is commensurable with 
$\prod_{i=1}^k\SL((n-1)h_i,\BZ(\zeta_{m_i}))^{m_i}$. 
\end{proposition}
\begin{proof}
The proof of Theorem \ref{teo4} shows that there is a subgroup of finite index
$\Gamma_0\le\A(F_n)$ and a representation 
\begin{equation}\label{Mech}
\rho_0 : \Gamma_0 \to \SL(n,\BZ)\times
\prod_{i=1}^k\SL((n-1)h_i,\BQ(\zeta_{m_i}))^{m_i} 
\end{equation}
such that the image $\rho_0(\Gamma_0)$ is the internal direct product of its
intersections with the factors in (\ref{Mech}). 
These intersections in turn have
finite index in the corresponding factor. Using the fact that a normal
subgroup of a subgroup of finite index in $\SL(n,\BZ)$ ($n\ge 3$) is either
finite or has finite index we conclude the proof of the prposition.
\end{proof}

Another easy consequence of Theorem \ref{teo4} concerns representations of the
the outer automorphism group ${\rm Out}(F_n)$. We have:
\begin{proposition}\label{teoo}
Let $n\ge 2$, $k\ge 1$, $h_1<\ldots <h_k$, $m_1,\ldots ,m_k$ be 
natural numbers.  Let $\BQ(\zeta_{m_i})$ be the field of 
$m_i$-th roots of unity and
$\BZ(\zeta_{m_i})$ its ring of integers. There is a 
subgroup $\Gamma\le {\rm Out}(F_n)$ of finite index and a representation 
$$
\rho : \Gamma \to \prod_{i=1}^k\SL((n-1)h_i,\BQ(\zeta_{m_i}))^{m_i} 
$$
such that $\rho(\Gamma)$ is commensurable with 
$\prod_{i=1}^k\SL((n-1)h_i,\BZ(\zeta_{m_i}))^{m_i}$. 
\end{proposition}
Here we only have to note that the image of the inner automorphisms under
any of the representations from Theorem \ref{teo4} is finite. This finite
subgroup of the image can be avoided by going to a subgroup of finite index in
$\Gamma$. 

As a consequence of Propositions \ref{teo3}, \ref{teoo} we get:
\begin{corollary}\label{coro8}
The groups ${\rm IA}(F_3)$ and ${\rm Out}(F_3)$ are large.
\end{corollary}

Our result from Theorem \ref{teo4} can 
conviniently be summarised in the language 
of the pro-algebraic completion ${\cal A}(\A(F_n))$ of $\A(F_n)$. 
For a definition and
discussion of properties of the pro-algebraic completion (also called the
Hochschild-Mostow group) of a group see
\cite{LM} and \cite {BLMM}. We have 
\begin{proposition}\label{teo9} 
Let $n\ge 2$ be a natural number and  ${\cal S}(\A(F_n))$ be
the semisimple part of the connected component of the identity of the 
pro-algebraic completion ${\cal A}(\A(F_n))$. Then for every $h\in \BN$, the
group $\SL((n-1)h,\BC)$ appears infinitely many times as a factor 
in ${\cal S}(\A(F_n))$. That is, there is a surjective homomorphism
$$
{\cal S}(\A(F_n))\to \prod_{h=1}^\infty \prod_{i=1}^\infty \SL((n-1)h,\BC).
$$
\end{proposition}
Propositions \ref{teo3}, \ref{teoo} imply similar results for the groups
${\rm IA}(F_n)$, ${\rm Out}(F_n)$ ($n\ge 3$).

\subsection{Representations of $\Gamma(G,\pi)$ into affine groups}\label{Lini}

In this section we give another construction of a linear representation
of subgroups of finite index in $\A(F_n)$. This representation takes values in 
affine groups (which have a non-trivial unipotent radical). The construction
originates from the proof of Gasch\"utz' result as presented in \cite{JR}.
Besides giving new types of image groups, it is of importance for the
computer algorithms which we have used to create the computational results
described in the next subsection.  

Let $G$ be a finite group and $\pi : F_n\to G$ a surjective homomorphism
of the free group $F_n$ onto $G$. Let $R$ be the kernel of 
$\pi$ and $\bar R$ the corresponding relation module. Let further
$\Gamma(G,\pi)$ be the subgroup of $\A(F_n)$ defined in 
(\ref{I1}). Given a group $H$ and a commutative ring $S$ we write $S[H]$ for
the corresponding group ring and ${\cal I}(S[H])$ for its augmentation ideal.
If $H_0\le H$ is a subgroup we define ${\cal I}(S[H],H_0)$ to be the two-sided
ideal of $S[H]$ generated by the $h-1$ for $h\in H_0$.

Gasch\"utz' result (\ref{gasch}) is proved (in \cite{JR}) 
by considering the exact sequence
\begin{equation}\label{gashpr}
0\to \BQ\otimes_\BZ \bar R\to {\cal I}(\BQ[F_n])/({\cal I}(\BQ[F_n],R)
\cdot{\cal I}(\BQ[F_n]))
\to {\cal I}(\BQ[G])\to 0
\end{equation}
of $G$-equivariant homomorphisms.
The right hand map is induced by $\pi : F_n\to G$ while the left hand map
comes from the map from $R$ to ${\cal I}(\BQ[F_n])$ which sends $r\in R$ to
$r-1$. The free group $F_n$ acts on the middle term in (\ref{gashpr}) by
multiplication from the left. This leads to an action of the 
finite group $G$ on this term. 

We now let $\Gamma(G,\pi)$ act on the left hand term as before, 
on the middle term
by its action on $F_n$ and trivially on $\BQ[G]$.
It is straightforward to see that the sequence
(\ref{gashpr}) is then $\Gamma(G,\pi)$-equivariant. We obtain a
representation
\begin{equation}\label{gashprep}
\eta_{G,\pi}:\Gamma(G,\pi)\to 
{\rm Hom}_G({\cal I}(\BQ[G]),\BQ\otimes_{\bb Z} \bar R)\rfish {\cal
  G}_{G,\pi}(\BQ).
\end{equation}
The semi-direct product in (\ref{gashprep}) is formed with respect to the
action of ${\cal G}_{G,\pi}(\BQ)$ on $\BQ\otimes_\BZ \bar R$.

Our methods show:
\begin{theorem}\label{teo16} Assume $n$ is a natural number with $n\ge 4$.
Let $\pi : F_n\to G$ be a redundant presentation of the finite group $G$. Then
$$\eta_{G,\pi}(\Gamma(G,\pi))\cap 
{\rm Hom}_G({\cal I}(\BQ[G]),\BQ\otimes_\BZ \bar R)\rfish 
{\cal G}^1_{G,\pi}(\BQ)$$ is 
of finite index in the arithmetic group 
${\rm Hom}_G({\cal I}(\BZ[G]),\bar R)\rfish{\cal G}_{G,\pi}^1(\BZ)$.
\end{theorem}
Theorem \ref{teo16} can then be used to prove
\begin{theorem}\label{teo17}
Let $n\ge 2$, $h$, $m$ be 
natural numbers.  Let $\BQ(\zeta_{m})$ be the field of 
$m$-th roots of unity and
$\BZ(\zeta_{m})$ its ring of integers. There is a 
subgroup $\Gamma\le \A(F_n)$ of finite index and a representation 
$$
\eta : \Gamma \to \left(\BQ(\zeta_m)^{(n-1)h}\rfish
\SL((n-1)h,\BQ(\zeta_{m}))\right)^{m} 
$$
such that $\eta(\Gamma)$ is commensurable with 
$(\BZ(\zeta_m)^{(n-1)h}\rfish  \SL((n-1)h,\BZ(\zeta_{m})))^{m}$. 
\end{theorem}

\subsection{Computational results}\label{Comp}

An important feature of our methods is that subgroups of quite high indices in
$\A(F_n)$ may be computationally located. We describe here some of the results
of our computer calculations. The computations where done using the computer
algebra system MAGMA. 

We shall describe our approach now in detail for $n=2$. Let 
$F_2:=\langle x,\, y\rangle$ be the free group. The automorphism group  
$\A^+(F_2)$ is generated by $\alpha,\, \beta$ which are given by
$$\alpha(x)=y^{-1},\quad \alpha(y)=x\qquad \beta(x)=x^{-1}y^{-1},\quad  
\beta(y)=x.$$
The group $\A^+(F_2)$ is finitely presented, a presentation is (see
\cite{Neu1}) 
\begin{equation}\label{neuma}
\A^+(F_2)=\langle\, \alpha,\, \beta \Mid \alpha^4=\beta^3=
\alpha^2\beta^2\alpha^2\beta\alpha\beta\alpha^2\beta^2\alpha=1\,\rangle.
\end{equation}
Let $\pi: F_2\to G$ be a presentation of a finite group. If $\varphi\in 
\A^+(F_2)$ is given by its values on $x,\, y$ it is easy to check whether
$\varphi$ is in $\Gamma(G,\pi)$ or not. Using a random process we generate
a list of words $\varphi_1,\ldots, \varphi_k\in \Gamma(G,\pi)$ in the
automorphisms $\alpha,\, \beta$. We then wait until
$$\Delta:=\langle \, \varphi_1,\ldots, \varphi_k\,\rangle$$
has finite index in $\A^+(F_2)$. This is checked using a Todd-Coxeter algorithm
using the presentation (\ref{neuma}). By this approach we first of all find the
following table. The first column contains the finite group $G$, which
we always think of being given as a permutation group. The third and fourth
column contain the images $x,\, y$ under the 
surjective homomorphism $\pi:F_2\to G$.
The fifth column gives the index of 
$\Delta\le \Gamma(G,\pi)$ in $\A^+(F_2)$. We give
the abelianised group $\Delta^{\rm ab}$ in the last column.   
\bigskip
\begin{center}
\begin{tabular}{|c|c|c|c|c|c|}
\hline
$G$  & $|G|$ & $\pi(x)$ & $\pi(y)$ & $[\A^+(F_2):\Delta]$ & $\Delta^{\rm ab}$\\
\hline
$C_2$ & 2 & (1,2) & (1) & 3 & $\BZ^2\times C_2\times C_4$\\
$C_3$ & 3 & (1,2,3) & (1) & 8 & $\BZ\times C_3^2$\\
$C_4$ & 4 & (1,2,3,4) & (1) & 12 & $\BZ^2\times C_4$\\
$C_2\times C_2$ & 4 & (1,2) & (3,4) & 6 & $\BZ^2\times C_2^3$\\
$C_5$ & 5 & (1,2,3,4,5) & (1) & 24 & $\BZ^3\times C_5$\\
$C_6$ & 6 & (1,2,3,4,5,6) & (1) & 24 & $\BZ^3\times C_6$\\
$S_3$ & 6 & (1,2,3) & (1,2) & 18 & $\BZ^2\times C_2$\\
$C_7$ & 7 & (1,2,3,4,5,6,7) & (1) & 48 & $\BZ^5\times C_7$\\
$D_4$ & 8 & (1,2,3,4) & (1,4)(2,3) & 24 & $\BZ^3\times C_2$\\
$Q_8$ & 8 & (1,7,2,8)(3,6,4,5) & (1,4,2,3)(5,7,6,8) & 24 & $\BZ^2\times C_4$\\
$D_5$ & 10 & (1,2,3,4,5) & (1,5)(2,4) & 30 & $\BZ^2\times C_2$\\
$A_4$ & 12 & (1,2,3) & (1,2)(3,4) & 96 & $\BZ^3$\\
$S_3\times C_2$ & 12 & (1,3)(4,5) & (1,2) & 36 & $\BZ^3\times C_2$\\
Sm($12,1$) & 12 & $\sigma_1$ & $\sigma_2$ & 72 & $\BZ^3\times C_2$\\
$A_5$ & 60 & (1,2,3,4,5) & (1,2,3) & 1080 & $\BZ^{17}$\\
\hline
\end{tabular}
\end{center}
The notation for the finite groups is: $C_n$ is the cyclic group of order $n$,
$D_n$ is the dihedral group of order $2n$,
$Q_8$ is the quaternion group of order $8$
and Sm($12,1$) is the non-abelian group of order $12$ not isomorphic to
$A_4$ or $S_3\times C_2$. The two permutations $\sigma_1,\, \sigma_2$ are:
$$\sigma_1:=(1, 8, 4, 11)(2, 9, 5, 12)(3, 7, 6, 10),\qquad
    \sigma_2:=(1, 9, 4, 12)(2, 7, 5, 10)(3, 8, 6, 11).$$

Proceeding with our computer calculations, we then evaluated 
the $\varphi_1,\ldots, \varphi_k$ (which now generate a
subgroup of finite index in $\Gamma(G,\pi)$) on the middle term of the sequence
(\ref{gashpr}). Locating the isotypic component of an irreducible
$\BQ[G]$-module $N$ inside all three terms of (\ref{gashpr}) we were able
to compute $\rho_{G,\pi,N}(\varphi_i)$ ($i=1,\ldots,k$) as matrices. Computing
determinants we could decide whether 
$\rho_{G,\pi}(\Gamma(G,\pi))\cap {\cal G}_{G,\pi}^1(\BZ)$ is of finite index
in $\rho_{G,\pi}(\Gamma(G,\pi))$ or not. 

We have considered at least one presentation $\pi: F_2\to G$ for
every group $G$ with $|G|\le 60$ which can be generated by 2 elements and
have found that $\rho_{G,\pi}(\Gamma(G,\pi))\cap {\cal G}_{G,\pi}^1(\BZ)$ is
of finite index in $\rho_{G,\pi}(\Gamma(G,\pi))$.   
Of course we have run experiments on many more groups 
(like $\PSL(2,7)$) but always found a
similar result. 

Let us report some more on the particularly interesting case $G=A_5$.
Let $N$ be the irreducible $\BQ[G]$-module of dimension 6. We have
${\rm End}_G(N)=\BQ(\sqrt{5})$. Its ring of integers is 
${\cal O}=\BZ((1+\sqrt{5})/2)$. 
Identifying the isotypic component of $N$ in $\BQ[G]$ with $M(3,\BQ(\sqrt{5}))$
we obtain a representation
$$\rho: \Gamma(G,\pi)\to \GL(1,M(3,\BQ(\sqrt{5})))=\GL(3,\BQ(\sqrt{5})).$$
Running the programm  described above we have found
\begin{proposition}\label{compu18}
For every presentation $\pi: F_2\to G=A_5$ there is a subgroup 
$\Delta\le \Gamma(G,\pi)$
of finite index such that $\rho_{A_5,\pi}(\Delta)$ 
is contained up to change of bases in the subgroup of 
$\SL(3,{\cal O})$ consisting of elements which have $(1,0,0)$ as a first row. 
\end{proposition}
Proposition \ref{compu18} shows that Theorem \ref{teo} is not true for $n=2$
and 
non-redundant presentations. It is not clear what Proposition \ref{compu18}
suggests toward the general case. It may suggest that Theorem \ref{teo} is not
true for non-redundant 
presentations but it may also be that $n=2$ is exceptional.

We have developed similar programms also for the cases $n=3$ and $n=4$.
They lead to the following tables. The notation is the same as in the above
table for $n=2$. The presentations used for $n=4$ take the value $1$ on the
fourth generator of the free group. 

\medskip
\centerline{\it Subgroups in $\A(F_3)$}
\begin{center}
\begin{tabular}{|c|c|c|c|c|c|c|}
\hline
$G$  & $|G|$ & $\pi(x)$ & $\pi(y)$ & $\pi(z)$ & $[\A(F_3):\Delta]$ & 
$\Delta^{\rm ab}$ \\
\hline
$C_2$ & 2 & (1,2) & (1) & (1) & 7 & $C_2^3$\\
$C_3$ & 3 & (1,2,3) & (1) & (1) & 26 & $C_6$ \\
$S_3$ & 6 & (1,2,3) & (1,2) & (1) & 168 & $C_2^3$\\
$D_4$ & 8 & (1,3) & (1,4,3,2) & (1) & 336 & $C_2^8$ \\
$Q_8$ & 8 & (1,7,2,8)(3,6,4,5) & (1,4,2,3)(5,7,6,8) & (1) & 336 & $C_2^7$\\
$D_5$ & 10 & (1,2,3,4,5) & (1,5)(2,4) & (1) & 840 & $C_2^5$ \\
$A_4$ & 12 & (1,2,3) & (1,2)( 3,4) & (1) & 1560 & $C_6$ \\
$S_3\times C_2$ & 12 & (1,3)(4,5) & (1,2) & (1) & 1008 & $C_2^9$\\
Sm($12,1$) & 12 & $\sigma_1$ & $\sigma_2$ & (1) & 1344 & $C_2^2\times C_4$\\
$A_5$ & 60 & (1,2,3,4,5) & (1,2,3) & (1) & 200160 & ? \\
\hline
\end{tabular}
\end{center}

\medskip
\centerline{\it Subgroups in $\A(F_4)$}
\begin{center}
\begin{tabular}{|c|c|c|c|}
\hline
$G$  & $|G|$ & $[\A(F_4):\Delta]$ & $\Delta^{\rm ab}$ \\
\hline
$C_2$ & 2 & 15 & $C_2^2$\\
$C_3$ & 3 & 80 & $C_6$\\
$S_3$ & 6 & 80 & $C_6$ \\
$D_4$ & 8 & 3360 & $C_2^4$ \\
$Q_8$ & 8 & 840 & $C_2^6$ \\
$D_5$ & 10 &  930 & $C_2^3$ \\
$A_4$ & 12 & 1680 & $C_3\times C_6$ \\
$S_3\times C_2$ & 12 & 2730 & $C_2^4$ \\
Sm($12,1$) & 12 & 3120 & $C_2^2\times C_4$\\
$A_5$ & 60 & 213098 & ? \\
\hline
\end{tabular}
\end{center}
We have run experiments on many finite groups $G$ and
presentations $\pi: F_n\to G$, we have always found
that $\rho_{G,\pi}(\Gamma(G,\pi))\cap {\cal G}_{G,\pi}^1(\BZ)$ is
of finite index in $\rho_{G,\pi}(\Gamma(G,\pi))$.



\begin{thebibliography}{99}

\bibitem{Ba} S.\, Bachmuth,
{\it Induced automorphisms of free groups and free metabelian groups.}
Trans. Amer. Math. Soc. {\bf 122}, (1966), 1--17

\bibitem{BLMM} H. Bass, A. Lubotzky, A.R. Magid, S. Mozes,
{\it The proalgebraic completion of rigid groups.}
Geometriae Dedicata {\bf 95}, (2002), 19--58

\bibitem{Bi} J.\, Birman,
{\it Braids, links and mapping class groups. \/},
Ann. of Math. Stud. {\bf 82}, Princeton Univ. Press, Princeton (1974)

\bibitem{CR} C.W.\, Curtis, I. Reiner,
{\it Representation Theory of finite Groups and associative Algebras.}
Interscience Publishers (1962) 

\bibitem{D} M. Deuring,
{\it Algebren.}
Ergebnisse {\bf 41}, Springer Verlag, New York, Berlin, Heidelberg (1968)

\bibitem{DJ} D. Djokovi\'c, V. Platonov,
{\it Low dimensional representations of $\A(F_2)$.}
Manuscripta Math. {\bf 89}, (1996), 475--509

\bibitem{DU} M. Dunwoody,
{\it Nielsen Transformations}. 
In: Computational Problems in Abstract Algebra (1970), Pergamon, Oxford,
45--46

\bibitem{EGM} J. Elstrodt, F. Grunewald, J. Mennicke,
{\it Groups acting on Hyperbolic Space: Harmonic Analysis
and Number Theory.}
Springer Monographs in Math., (1998), 524 pp 

\bibitem{Ev} M. Evans,
{\it Presentations of groups involving more generators than are necessary.}
Proc London Math. Soc. {\bf 67}, (1993), 106--126

\bibitem{Ev1} M. Evans,
{\it T-systems of certain finite simple groups.}
Math. Proc. Cambridge Philos. Soc. {\bf 113}, (1993), 9--22

\bibitem{F} E. Formanek,
{\it Characterizing a free group in its automorphism group.}
J. Algebra, {\bf 133}, (1990), 424--432

\bibitem{FP} E. Formanek, C. Procesi,
{\it The automorphism group of a free group is not linear.}
J. Algebra, {\bf 149}, (1992), 494--499

\bibitem{G} W. Gasch\"utz,
{\it \"Uber modulare Darstellungen endlicher Gruppen, die von freien Gruppen
  induziert werden.}
Math. Z. {\bf 60}, (1954), 274--286

\bibitem{Gi} R. Gilman,
{\it Finite quotients of the automorphism group of a free group.}
Canad. J. Math. {\bf 29}, (1977), 541--551

\bibitem{Gr} K. Gruenberg,
{\it Relation modules of finite groups.}
CBMS {\bf 25}, Amer. Math. Soc. (1976)

\bibitem{GLu} F. Grunewald,, A. Lubotzky,
{\it Linear representations of the mapping class group.}
In preparation

\bibitem{GP} F. Grunewald, V. Platonov,
{\it Rigidity results for groups with radical, cohomology of finite groups and
  arithmeticity problems.}
Duke Math. J. {\bf 100}, (1999), 321--358

\bibitem{GS} R. Guralnick, A. Shalev,
{\it On the spread of finite simple groups.}
Combinatorica {\bf 23}, (2003), 73--87 

\bibitem{JR} M. Jarden, J. Ritter, 
{\it Normal Automorphisms of absolute Galois Groups of $\wp$-adic Fields.}
Duke Math. J. {\bf 47}, (1980), 47--56
 
\bibitem{L} A. Lubotzky,
{\it Torsion in profinite completions of torsion-free groups.}
Quart. J. Math. Oxford Series {\bf 44}, (1993), 327--332

\bibitem{LM} A. Lubotzky, A.R. Magid, 
{\it Varieties of representations of finitely generated groups.}
Mem. Amer. Math. Soc. {\bf 58}, (1985)

\bibitem{LP} A. Lubotzky, I. Pak,
{\it The product replacement algorithm and Kazhdan's property (T).}
J. Amer. Math. Soc. {\bf 14}, (2001), 347--363 

\bibitem{MKS} W.\, Magnus, A. Karass, D. Solitar,
{Combinatorial Group Theory.}
Interscience Publishers (1966)

\bibitem{Neu1} B.H. Neumann,
{\it Die Automorphismengruppe der freien Gruppen.}
Math. Ann. {\bf 107}, (1933), 367--386

\bibitem{Neu} B.H. Neumann,
{\it On a question of Gasch\"utz.}
Archiv der Math. {\bf 7}, (1956), 87--90

\bibitem{N} J. Nielsen,
{\it Die Isomorphismengruppe der freien Gruppen.}
Math. Ann. {\bf 91}, (1924), 169--208

\bibitem{Pak} I. Pak,
{\it What do we know about the product replacement algorithm?}
Groups and Computation III (Columbus, Ohio, 1999), Ohio State
Univ. Math. Res. Inst. Publ. {\bf 8}, 301--347, de Gruyter, Berlin, (2001)

\bibitem{PLR} V. Platonov, A. Rapinchuk,
{\it Algebraic Groups and Number Theory.} Pure and Applied Mathematics {\bf
  139}, Academic Press, (1994)

\bibitem{PORA} A. Potapchik, A. Rapinchuk,
{\it Low-dimensional linear representations of $\A(F_n)$, $n\ge 3$.}
Transactions of the Amer. Math. Soc. {\bf 352}, (1999), 1437--1451 

\bibitem{RA} A. Rapinchuk,
{\it $n$-dimensional linear representations of $\A(F_n)$, and more.}
J. of Algebra {\bf 197}, (1997), 146--152

\bibitem{Serresl} J.-P., Serre,
{\it Le probl\'eme des groupes des congruence pour $\SL_2$.}
Ann. of Math., {\bf 92}, (1970), 489--527

\bibitem{S} J.-P., Serre,
{\it Linear Representations of Finite Groups.}
Graduate Text in Math. {\bf 42},
Springer Verlag, New York, Berlin, Heidelberg (1993)

\bibitem{Vas} L.\, Vaserstein,
{\it The structure of classical arithmetic groups of rank greater than one.}
Math. USSR Sbornik (English translation), {\bf 20}, (1973), 465--492 

\bibitem{Vas1} L.\, Vaserstein,
{\it The group $\SL_2$ over Dedekind rings of arithmetic type.}
Mat. Sb. {\bf 89}, (1972), 313--322

\bibitem{Venk} T.N.\, Venkataramana,
{\it On systems of generators of arithmetic subgroups of higher rank groups.}
Pacific Journal of Math., {\bf 166, 1}, (1994), 193--212

\end{thebibliography}
\end{document}